\documentclass[preprint,12pt]{elsarticle}
\usepackage[letterpaper,margin=1in]{geometry}

\usepackage{graphicx}
\usepackage[export]{adjustbox}
\usepackage{caption}
\usepackage{subcaption}
\usepackage{amssymb,amsmath,amsthm,amsfonts,amscd, bm}
\usepackage{hyperref}
\usepackage{mathtools}
\usepackage{algorithm}
\usepackage{algorithmicx}
\usepackage{algpseudocode}
\usepackage{xcolor}
\usepackage{booktabs}
\usepackage{array}
\newcolumntype{C}[1]{>{\centering\arraybackslash}p{#1}}
\newcolumntype{P}[1]{>{\raggedright\arraybackslash}p{#1}}
\usepackage{ragged2e}
\usepackage{url}
\usepackage{multirow}
\usepackage{makecell}
\usepackage{tabularx}
\usepackage{xparse, etoolbox}
\usepackage{tikz}
\usetikzlibrary{arrows.meta,positioning,calc}
\hypersetup{colorlinks=false}
\newtheorem{thm}{Theorem}[section]

\newtheorem{lem}[thm]{Lemma}
\newtheorem{prop}[thm]{Proposition}

\newtheorem{assumption}[thm]{Assumption}

\theoremstyle{definition}

\theoremstyle{remark}

\newcommand\lr[1]{\left({#1}\right)}
\newcommand\norm[1]{\left\lVert#1\right\rVert}
\newcommand{\abs}[1]{\left\lvert#1\right\rvert}

\DeclareMathOperator*{\argmin}{arg\,min}
\DeclarePairedDelimiterX\innerp[2]{\langle}{\rangle}{#1,#2}
\DeclarePairedDelimiterX{\Innerp}[1]{\langle}{\rangle}{\Innpargs{#1}}

\begin{document}
\makeatletter
\def\ps@pprintTitle{%
   \let\@oddhead\@empty
   \let\@evenhead\@empty
   \let\@oddfoot\@empty
   \let\@evenfoot\@oddfoot}
\makeatother

\begin{frontmatter}

\title{Learning Hamiltonian Flow Maps from Numerical-Scheme Residuals for Long-Time Multiscale Simulation}

\author[1]{Rui Fang and Richard Tsai}
\affiliation[1]{organization={Department of Mathematics and Oden Institute for Computational Engineering and Sciences},
            addressline={The University of Texas at Austin}, 
            city={Austin},
            postcode={78712}, 
            state={TX},
            country={USA}}

\begin{abstract}

Hamiltonian systems with widely separated timescales arise in molecular dynamics, classical mechanics, and plasma physics. Long-time simulation of such systems is expensive because standard direct integrators generally need to resolve the fastest dynamics even when the quantities of interest evolve on much slower scales. The cost is particularly severe for large ensembles of trajectories.

We develop a framework for learning Hamiltonian flow maps directly with neural networks. The map is trained either from the residual of a convergent numerical scheme or from reference trajectory data. For variable-timestep maps, truncated Taylor expansions enforce the correct short-time behavior while a neural network represents the remainder. For stiff oscillatory systems, the losses are measured in an energy-balanced norm that weights position errors according to their associated frequencies. We also introduce an HMC-$H_0$ procedure for generating training samples from microcanonical energy surfaces.

The analysis clarifies what the residual training learns and how local flow-map errors propagate. Under a nonsingularity condition on the implicit part of a one-step scheme, every $C^1$ critical point of the unrestricted function-space residual functional has zero residual and therefore reproduces that numerical scheme. For a class of multiscale oscillatory Hamiltonians, the amplification rate of recursively applied flow-map errors is independent of the stiff frequencies when the local error is controlled in the energy-balanced norm; linear accumulation follows on the corresponding pre-asymptotic time window. Numerical experiments on separable, nonseparable, and noncanonical systems demonstrate long-time accuracy and identify regimes in which learned flow maps can reduce the cost of multiscale ensemble simulation.

\end{abstract}



\begin{keyword}
Hamiltonian systems \sep flow maps \sep geometric numerical integration \sep long-time simulation \sep multiscale dynamics \sep deep learning


\end{keyword}

\end{frontmatter}

\section{Introduction}

Hamiltonian systems are governed by ordinary differential equations of the form
\begin{equation}\label{eq:generic-Hamiltonian-system}
  \frac{d}{dt}u=f(u)\equiv J^{-1}\nabla H(u),\qquad
  J=\begin{pmatrix}0&I\\-I&0\end{pmatrix},\qquad
  u=(p,q)\in\mathbb{R}^{2d},
\end{equation}
where $p$ and $q$ denote generalized momenta and positions and $H$ is the Hamiltonian. The exact flow preserves the symplectic form, phase-space volume, and the energy. Symplectic numerical integrators exploit this geometry and often provide substantially better long-time behavior than general-purpose methods; see, e.g., \cite{feng1986difference,sanz1992symplectic,geometric}.

The computational difficulty increases sharply in multiscale systems. If $H=H_\epsilon$ contains fast dynamics on a timescale $\epsilon\ll1$, a direct integrator must generally resolve those dynamics even when the quantities of interest evolve on an $\mathcal{O}(1)$ or slower scale. Over long horizons, and especially for large ensembles of initial conditions, this requirement can dominate the computational cost. Classical approaches for reducing it include heterogeneous multiscale methods \cite{engquist2005heterogeneous,engquist2007heterogeneous,weinan2003heterognous,ariel2009multiscale}, stroboscopic averaging \cite{calvo2011stroboscopic,calvo2011numerical,chartier2016solving}, and time-parallel multiscale methods \cite{ariel2016parareal,ariel2017theta,nguyen2020stable}.

A complementary strategy is to learn the flow map itself. A fixed-timestep model approximates $u\mapsto\phi_{T_0}(u)$ for a prescribed increment $T_0$, whereas a variable-timestep model approximates $(u,t)\mapsto\phi_t(u)$. Once trained, such a model can advance an ensemble by direct network evaluation, without numerically integrating the vector field at inference time. Many existing approaches impose symplecticity through the architecture. Representative examples include SympNets \cite{jin2020sympnets}, H\'enonNets \cite{burby2020fast}, generating-function networks \cite{chen2021data}, SympFlow \cite{canizares2024symplectic}, and the SymplecticGyroceptron \cite{duruisseaux2023approximation}; see \cite{celledoni2021structure} for a broader account of structure-preserving deep learning. These methods provide powerful architectural mechanisms for geometric structure. Our approach is complementary: we use generic feedforward networks and place the numerical structure primarily in the loss, the short-time parameterization, and the sampling distribution.

\paragraph{Our approach}
We learn either a fixed-timestep map $\Phi_{T_0}(u)$ or a variable-timestep map $\Phi(u,t)$. The latter is constructed from a truncated Taylor expansion at $t=0$ plus a learned remainder. Training can be unsupervised, by minimizing the residual of a chosen convergent one-step numerical scheme, or supervised, using reference trajectory data computed offline. For stiff systems, both losses can be measured in an energy-balanced norm that resolves errors in fast positions on their natural energetic scale. Training states can be sampled from microcanonical energy surfaces with the HMC-$H_0$ procedure introduced in Section~\ref{sec:method_data}.

The main contributions are:
\begin{enumerate}
    \item \emph{Scheme-based residual losses and their critical points.} We show that, for a broad class of one-step schemes, every $C^1$ critical point of the unrestricted function-space residual functional has zero residual whenever the Jacobian with respect to the new state is nonsingular. This condition is automatic for explicit methods and for the Velocity Verlet splitting used here; for genuinely implicit schemes it is a solvability condition on the implicit update, not a stability condition (Section~\ref{sec:critical-pts-of-R}).
    \item \emph{Taylor-consistent neural flow maps.} Generic feedforward networks are parameterized so that the variable-timestep map has prescribed consistency order at $t=0$, while allowing different Taylor orders for slow and fast variables (Section~\ref{sec:method_arch}).
    \item \emph{Energy-balanced losses and invariant-measure sampling.} We introduce frequency-weighted losses for stiff oscillatory Hamiltonians and HMC-$H_0$ sampling on prescribed energy surfaces. For the multiscale class analyzed in Section~\ref{sec:ana_errorbound}, the weighted norm yields an error-amplification rate independent of the stiff frequencies; linear accumulation follows when the corresponding slow-scale amplification remains small.
    \item \emph{Numerical assessment across multiscale systems.} We test the framework on nearly periodic coupled oscillators, the Fermi--Pasta--Ulam--Tsingou (FPUT) problem, and a noncanonical $\alpha$-particle model. The experiments include ablations, a parameter-matched symplectic-network comparison, and runtime measurements that identify when learned flow maps are computationally advantageous (Section~\ref{sec:dl-results}).
\end{enumerate}

Figure~\ref{fig:overview} summarizes the relation among the training distribution, loss, architecture, theory, and numerical validation.

\begin{figure}[t!]
\centering
\resizebox{\textwidth}{!}{%
\begin{tikzpicture}[
  font=\small,
  >={Stealth[length=2.4mm]},
  pillar/.style={draw=black!75, thick, rounded corners=2.5pt, fill=black!4,
    text width=4.6cm, inner sep=7pt, align=flush left, minimum height=3.45cm},
  theory/.style={draw=black!45, rounded corners=2pt, fill=blue!8,
    text width=4.6cm, inner sep=5pt, align=flush left, font=\footnotesize},
  stage/.style={draw=black!75, thick, rounded corners=2.5pt, fill=black!10,
    inner sep=7pt, align=center},
  flow/.style={->, thick, black!70}
]
\node[pillar, anchor=north] (loss) at (0,0) {%
  {\footnotesize\itshape what to fit, in which norm}\\[1pt]
  \textbf{Loss functions \; (Sec.~\ref{sec:dl})}\\[3pt]
  scheme-based residual $R_h[\Phi]$, or data loss on examples $(u,\phi_{T_0}(u))$;\\[1pt]
  both measured in the energy-balanced norm $\|\cdot\|_W$ (Sec.~\ref{sec:method_energy_norm})};
\node[pillar, anchor=north] (data) at (-5.55,0) {%
  {\footnotesize\itshape where in phase space to fit}\\[1pt]
  \textbf{Training data \; (Sec.~\ref{sec:method_data})}\\[3pt]
  HMC-$H_0$: microcanonical samples on target energy surfaces};
\node[pillar, anchor=north] (arch) at (5.55,0) {%
  {\footnotesize\itshape with which function class}\\[1pt]
  \textbf{Neural flow maps \; (Sec.~\ref{sec:method_arch})}\\[3pt]
  Taylor-consistent networks $\Phi = u + \Psi^{(p)}$; fixed-timestep $\Phi_{T_0}$, variable-timestep $\Phi(u,t)$; slow--fast orders $(p_s,p_f)$};
\draw[flow] ($(data.east)+(0,0.9)$) -- ($(loss.west)+(0,0.9)$);
\draw[flow] ($(arch.west)+(0,0.9)$) -- ($(loss.east)+(0,0.9)$);
\node[theory, below=3mm of data] (thdata) {%
  \textbf{Prop.~\ref{prop:hmc-invariance}}\; the sampler leaves the microcanonical measure invariant};
\node[theory, below=3mm of loss] (thloss) {%
  \textbf{Thm.~\ref{thm:general-one-step}}\; no nonzero-residual critical points in function space\\[1pt]
  \textbf{Thm.~\ref{thm:linear-error-growth}}\; stiffness-uniform error propagation in $\|\cdot\|_W$};
\node[theory, below=3mm of arch] (tharch) {%
  \textbf{Lem.~\ref{lem:consistency}}\; built-in consistency of order $p$ at $t=0$};
\node[stage, below=9mm of thloss, text width=5.6cm] (trained) {%
  \textbf{Trained flow map}\\[1pt]
  recursive application $\Phi^{(n)}$\\ for long-time simulation};
\node[theory, text width=4.2cm] (thtrained) at (trained -| arch) {%
  \textbf{Thm.~\ref{thm:linear-error-growth}}\; amplification rate independent of the stiff frequencies};
\draw[flow] (thdata.south) to[out=-85, in=180] (trained.west);
\draw[flow] (thloss.south) -- (trained.north);
\draw[flow] (tharch.south) to[out=-95, in=55] (trained.north east);
\draw[black!45, thick] (trained.east) -- (thtrained.west);
\node[stage, below=7mm of trained, text width=15.3cm, fill=black!14] (valid) {%
  \textbf{Validation on challenging multiscale benchmarks \; (Sec.~\ref{sec:dl-results})}\\[2pt]
  \mbox{NPCO (slow--fast) \;$\cdot$\; FPUT, $\omega=50$ and $300$ (stiff, chaotic) \;$\cdot$\; $\alpha$-particle (noncanonical)}\\[1pt]
  a parameter-matched symplectic-network baseline, scope of existing learned integrators, runtime (Secs.~\ref{sec:baselines}--\ref{subsec:runtime-comp})};
\draw[flow] (trained.south) -- (valid.north);
\end{tikzpicture}}%
\caption{Overview of the proposed framework. Each component of the learning pipeline---training data, loss functions, and network architecture---is paired with a corresponding theoretical guarantee (shaded tags), and the assembled pipeline is validated on multiscale benchmarks spanning substantial scale separation and noncanonical dynamics.}
\label{fig:overview}
\end{figure}

\paragraph{Why scheme-based residuals?}
A scheme residual couples the learned map at times separated by a nonzero increment $h$. It therefore enforces a discrete dynamical relation without differentiating the network with respect to time. In the variable-timestep experiments below, this makes residual evaluation cheaper than a pointwise PINN residual based on automatic differentiation and allows the training constraint to act over a finite temporal interval. The price is equally clear: at zero residual the network reproduces the chosen numerical scheme, not the exact flow. The discretization error of that scheme must therefore be controlled separately; Section~\ref{sec:critical-pts-of-R} makes this distinction precise. Related numerical-scheme least-squares formulations for Hamilton--Jacobi equations were studied in \cite{esteve2025finite,bokanowski2026solving}. The former treats finite-difference residual minimization for a Lax--Friedrichs-type discretization, while the latter develops broader well-posedness and critical-point results for monotone discretizations. Those works concern discrete Hamilton--Jacobi equations; here we analyze continuous-in-time one-step residuals for learned flow maps.

\section{Learning objectives and numerical-scheme residuals}\label{sec:dl}

Recall the generic Hamiltonian system in \eqref{eq:generic-Hamiltonian-system}.  
Given an initial condition $u_0 \in \Omega \subset \mathbb{R}^{2d}$ and $t \in [0,T]$, the exact flow map $\phi_t(u_0)$ returns the initial value problem's solution at time $t$.  
Our goal is to construct a neural network  
\[
    \Phi: \Omega \times [0,T] \to \mathbb{R}^{2d}, \quad \Phi(u,t) \approx \phi_t(u),
\]
satisfying the initial condition
\begin{equation}
    \Phi(u,0)=u,\qquad \forall u\in\Omega.
\end{equation}

We determine $\Phi$ by solving the variational problem  
\begin{align}
    \Phi^* = \argmin_{\Phi \in \mathcal{X}} \mathcal{R}[\Phi], 
    \label{eq:total_loss_functional}
\end{align}
where $\mathcal{X}$ is the function space defined by the network architecture (to be discussed in Section~\ref{sec:method_arch}), and
the residual loss
\begin{align}\label{eq:var_timestep_loss_functional}
    \mathcal{R}[\Phi] &:= \tfrac12 \int_{\Omega} \int  \|R_h[\Phi](u, t)\|^2 \, d\nu(t) \, d\rho(u)
\end{align}
encodes a notion of consistency with the governing differential equations, where $\nu$ is a probability measure on $[0,T]$ and $\rho$ is a probability measure on the phase-space domain $\Omega$. The analysis of Section~\ref{sec:critical-pts-of-R} takes $d\nu=\gamma(t)\,dt$ with $\gamma$ continuous and strictly positive. This covers the population objective associated with continuously distributed, and in particular uniformly resampled, collocation times. A fixed finite grid gives a discrete quadrature approximation and is not covered directly by the theorem. The residual $R_h[\Phi]$ is obtained by writing a numerical integrator as a one-step relation with step size $h$. Section~\ref{sec:critical-pts-of-R} analyzes the corresponding unrestricted function-space functional. 


We also consider the approximation of the flow map at a fixed time increment $T_0$. This will involve data of the form $(u, \phi_{T_0}(u))$ and the minimization of the data loss with respect to a probability measure $\rho_H$ over $\Omega$:
\begin{align}\label{def:data-loss}
    \mathcal{L}[\Phi_{T_0}; \phi] &:= \frac{1}{2S}   \sum_{k=1}^{S} \int_{\Omega} \norm{ \Phi^k_{T_0}(u)  - \phi_{kT_0}(u)}^2 \, d\rho_H(u), 
\end{align}
where 
\begin{align*}
    \Phi^k_{T_0} :=\underbrace{\Phi_{T_0}\circ \Phi_{T_0}\circ \cdots \circ\Phi_{T_0}}_{k\text{ times}},
\end{align*}
so that $\Phi^k_{T_0}(u)$ approximates the state at time $kT_0$ starting from $u$.

\subsection{Energy-balanced norms}\label{sec:method_energy_norm}

For stiff oscillatory systems, the Euclidean norm can underweight errors in fast position variables. Consider a separable oscillatory component with momenta $p$, positions $q$, and frequencies $\omega_j$. Define
\[
\tilde\omega_j:=\max\{\omega_j,1\},\qquad
\tilde\Omega:=\operatorname{diag}(\tilde\omega_1,\ldots,\tilde\omega_d),
\]
and the energy-balanced norm
\begin{equation}\label{def:weighted-l2-norm}
    \norm{(v_p,v_q)}_W
    :=\norm{(v_p,\tilde\Omega v_q)}_2
    =\norm{W(v_p,v_q)}_2,
    \qquad
    W:=\operatorname{diag}(I,\tilde\Omega).
\end{equation}
A displacement in a fast position $q_j$ is therefore weighted by its frequency. This is the natural scaling for the quadratic energy $\frac12\omega_j^2q_j^2$: on an $\mathcal O(1)$ energy surface, $q_j=\mathcal O(\omega_j^{-1})$.

Whenever this scaling is appropriate, we replace the Euclidean norms in the data and residual objectives by $\norm{\cdot}_W$. Thus
\begin{align}
    \mathcal L_W[\Phi_{T_0};\phi]
    &:=\frac{1}{2S}\sum_{k=1}^S\int_\Omega
      \norm{\Phi_{T_0}^k(u)-\phi_{kT_0}(u)}_W^2\,d\rho_H(u),
    \label{def:energy-balanced-data-loss}\\
    \mathcal R_W[\Phi]
    &:=\frac12\int_\Omega\int
      \norm{R_h[\Phi](u,t)}_W^2\,d\nu(t)\,d\rho(u).
    \label{def:energy-balanced-residual-loss}
\end{align}
For the FPUT variables introduced in Section~\ref{subsec:FPU}, ordered as $u=(y_s,x_s,y_f,x_f)$ with common fast frequency $\omega>1$, this becomes
\[
    W=\operatorname{diag}(I_{3m},\omega I_m),
\]
which is the weighting used in the experiments. Section~\ref{sec:ana_errorbound} explains why this norm removes the stiff frequencies from the error-amplification constant. The same principle can be generalized to other sensitive directions; nonlinear transformations of this type are considered in \cite{nguyen2020stable,nguyen2023numerical,fang2024stabilization}.

Because these maps are applied recursively, it is useful to concentrate $\rho$ and $\rho_H$ on dynamically relevant regions. Section~\ref{sec:method_data} develops this idea using microcanonical measures on energy surfaces.

\subsection{Residuals of numerical schemes}
\label{sec:method_var_timestep_training}

The residual measures the discrepancy between the learned flow map and the numerical update defined by a chosen scheme with step size $h$. For example, forward Euler has the update map $\Phi_h^{\mathrm{FE}}(u)=u+hf(u)$, giving
\begin{align}
    R_h[\Phi](u,t)
    :=\frac1h\lr{\Phi(u,t+h)-\Phi_h^{\mathrm{FE}}(\Phi(u,t))}.
    \label{eq:forward_euler}
\end{align}
A one-step residual couples times separated by $h$. When the time measure is Lebesgue and $T=Nh$, the continuous loss can be decomposed, up to an irrelevant normalization factor, into shifted grids. With $t_n=nh$, $n=0,\ldots,N-1$, define
\begin{equation}\label{eq:var_timestep_loss_functional_grid_shifted}
    \widehat{\mathcal R}_h[\Phi;\tau]
    :=\frac12\int_\Omega\sum_{n=0}^{N-1}
      \norm{R_h[\Phi](u,t_n+\tau)}^2\,d\rho(u),
    \qquad 0\le\tau<h.
\end{equation}
Then
\[
    \frac12\int_\Omega\int_0^T\norm{R_h[\Phi](u,t)}^2\,dt\,d\rho(u)
    =\int_0^h\widehat{\mathcal R}_h[\Phi;\tau]\,d\tau.
\]
For a general density $d\nu=\gamma(t)\,dt$, the same decomposition holds after inserting the factor $\gamma(t_n+\tau)$ inside the sum defining $\widehat{\mathcal R}_h$; the outer measure remains $d\tau$.
In computation, the time integral is replaced by finitely many collocation times. A single fixed grid is only a quadrature approximation to the continuous functional: even zero residual at the grid points need not imply small residual between them. The following example compares this behavior with a pointwise differential-equation residual.

\begin{figure}[ht]
    \centering
    \includegraphics[width=0.6\linewidth]{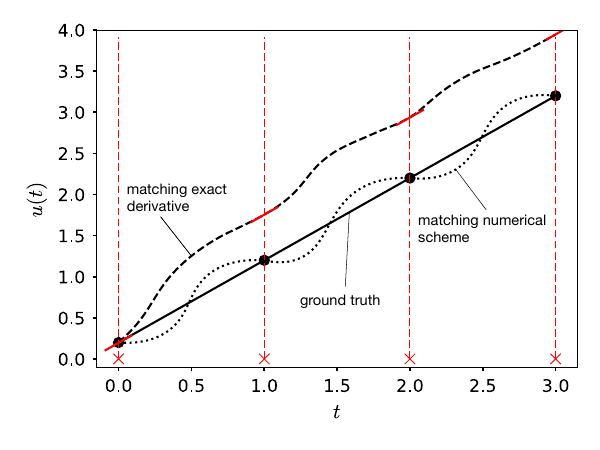}
    \caption{Two residual constructions at time-collocation points. The exact ODE residual is local in time, whereas a one-step scheme residual couples the network values at times separated by $h$ (red crosses). With finitely many collocation points, neither loss directly constrains every intermediate time; the distinction is that the scheme residual imposes a nonlocal relation between paired times. Section~\ref{sec:critical-pts-of-R} determines what the corresponding continuous residual functional enforces.}
    \label{fig:fd_illustration}
\end{figure}

\paragraph{A motivating example}  
We illustrate the method on the harmonic oscillator  
\[
H(p,q) = \tfrac12 p^2 + \tfrac12 q^2, 
\quad \dot{p} = -q, \quad \dot{q} = p,
\]
in a 2D phase space.  As the scheme-based residual, we use the Velocity Verlet (VV) integrator, a standard symplectic method for separable systems $H(p,q)=\tfrac{1}{2}p^{\!\top}M^{-1}p+V(q)$ with force $F(q)=-\nabla V(q)$. Its one-step update reads
\[
p_{n+\frac12}=p_n+\tfrac{h}{2}F(q_n),\qquad 
q_{n+1}=q_n+hM^{-1}p_{n+\frac12},\qquad 
p_{n+1}=p_{n+\frac12}+\tfrac{h}{2}F(q_{n+1}),
\]
giving the residual $R_h[\Phi](u,t) := \frac{1}{h}\lr{\Phi(u,t+h) - \Phi^{\text{VV}}_h(\Phi(u,t))}$.
We compare two residual definitions: (1) VV with $h=0.5$; (2) the exact (ODE) residual $R_{\text{exact}}[\Phi](u,t) := \partial_t \Phi(u,t) - J^{-1}\nabla H(\Phi(u,t))$.  

Each model is trained for 50{,}000 iterations using 10 phase-space points uniformly sampled on $H(p,q) = \tfrac12$ and $N$ time-collocation points in $[0,10]$, either on a fixed grid ($N=11,21,41$) or randomly resampled at each iteration ($N=41$).
With $N=41$, training takes approximately $394$\,s for VV residuals and $702$\,s for the exact residual, showing that automatic differentiation of the latter is more expensive.

\begin{figure}
    \centering
    \begin{tabular}{@{}c@{\hspace{2pt}}c@{\hspace{2pt}}c@{\hspace{2pt}}c@{\hspace{2pt}}c@{}}
     & \multicolumn{2}{c}{VV residual, $h=0.5$} & \multicolumn{2}{c}{exact residual}\\
        \rotatebox{90}{\small\hspace{1.6em} $N=11$} &
        \begin{subfigure}[b]{0.215\linewidth}
            \centering
            \includegraphics[width=\linewidth,trim=0 10 180 10, clip]{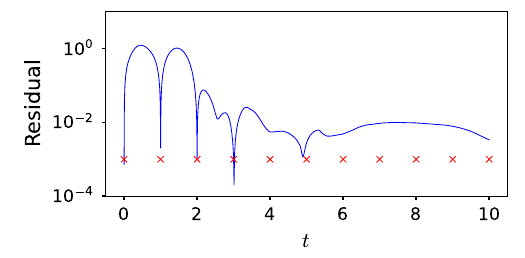}
        \end{subfigure} &
        \begin{subfigure}[b]{0.215\linewidth}
            \centering
            \includegraphics[width=\linewidth,trim=90 10 190 10, clip]{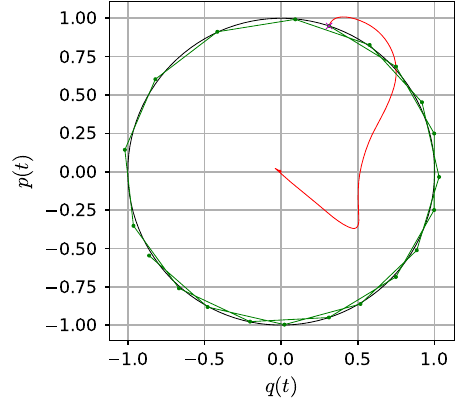}
        \end{subfigure} &
        \begin{subfigure}[b]{0.215\linewidth}
            \centering
            \includegraphics[width=\linewidth,trim=0 10 180 10, clip]{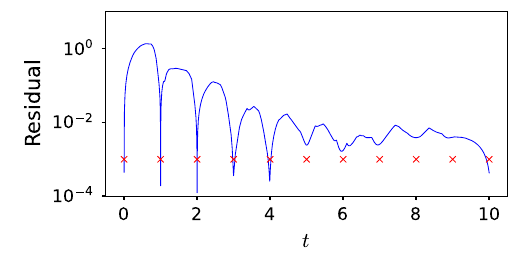}
        \end{subfigure} &
        \begin{subfigure}[b]{0.215\linewidth}
            \centering
            \includegraphics[width=\linewidth,trim=100 10 190 10, clip]{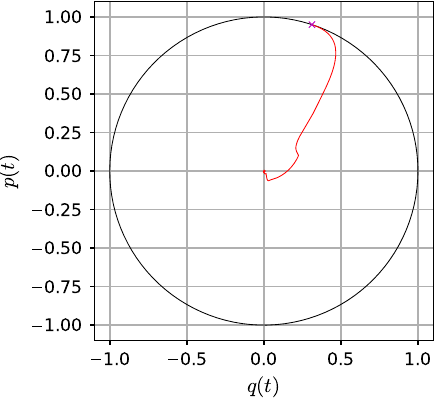}
        \end{subfigure} \\
        \rotatebox{90}{\small\hspace{1.6em} $N=21$} &
        \begin{subfigure}[b]{0.215\linewidth}
            \centering
            \includegraphics[width=\linewidth,trim=0 10 180 10, clip]{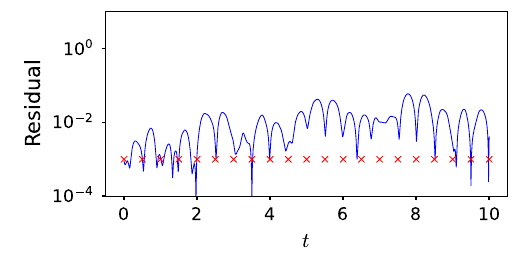}
        \end{subfigure} &
        \begin{subfigure}[b]{0.215\linewidth}
            \centering
            \includegraphics[width=\linewidth,trim=90 10 190 10, clip]{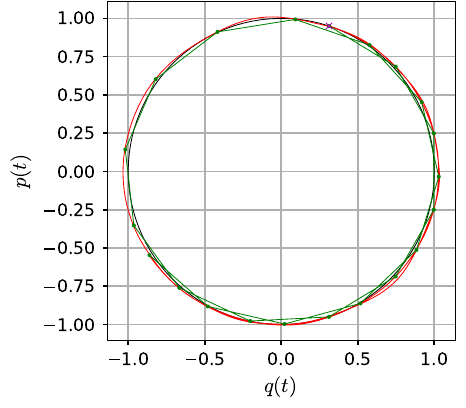}
        \end{subfigure} &
        \begin{subfigure}[b]{0.215\linewidth}
            \centering
            \includegraphics[width=\linewidth,trim=0 10 180 10, clip]{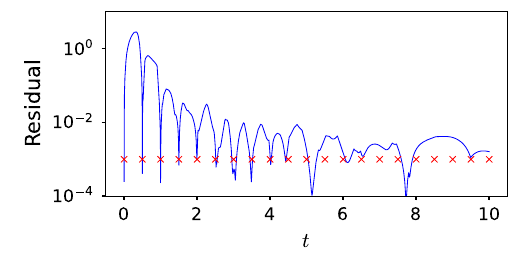}
        \end{subfigure} &
        \begin{subfigure}[b]{0.215\linewidth}
            \centering
            \includegraphics[width=\linewidth,trim=100 10 190 10, clip]{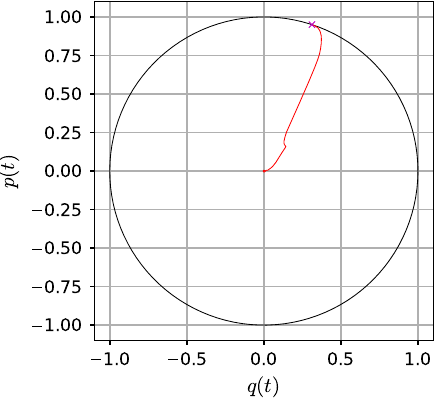}
        \end{subfigure} \\
        \rotatebox{90}{\small\hspace{1.6em} $N=41$} &
        \begin{subfigure}[b]{0.215\linewidth}
            \centering
            \includegraphics[width=\linewidth,trim=0 10 180 10, clip]{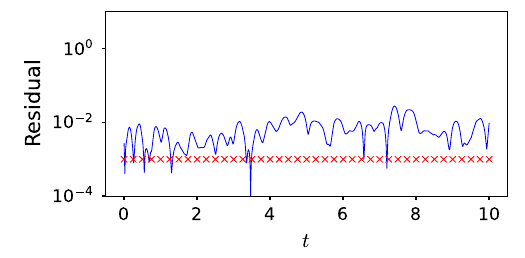}
        \end{subfigure} &
        \begin{subfigure}[b]{0.215\linewidth}
            \centering
            \includegraphics[width=\linewidth,trim=90 10 190 10, clip]{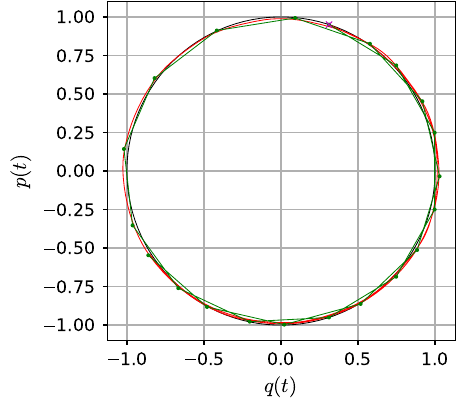}
        \end{subfigure} &
        \begin{subfigure}[b]{0.215\linewidth}
            \centering
            \includegraphics[width=\linewidth,trim=0 10 180 10, clip]{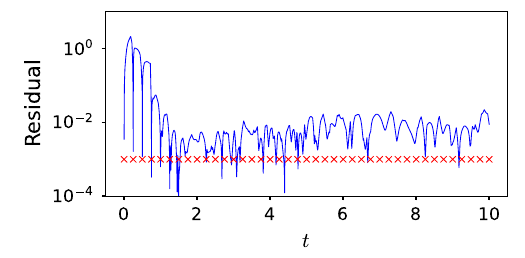}
        \end{subfigure} &
        \begin{subfigure}[b]{0.215\linewidth}
            \centering
            \includegraphics[width=\linewidth,trim=100 10 190 10, clip]{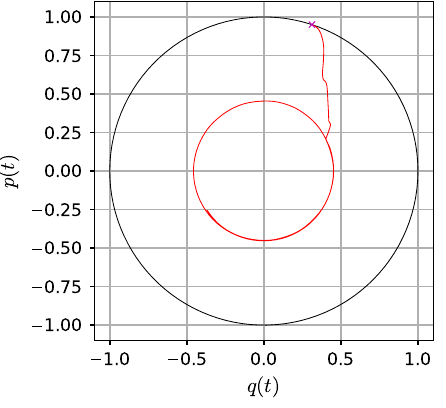}
        \end{subfigure}
    \end{tabular}
    \caption{(Harmonic oscillator) Residual-vs-time plots and phase-space trajectories of one-step predictions generated by flow maps trained with the velocity Verlet residual at $h=0.5$ (left pair) and with the exact residual (right pair), using $N=11$, $21$ or $41$ evenly spaced time-collocation points on $[0,10]$. In the residual plots, blue curves indicate the residual values and red crosses the collocation points. In the phase-space plots, the black curve is the exact flow, the red curve the learned flow, and the green points the numerical solution at discrete times ($t=0,h,2h,\dots$). Away from the collocation points the exact-residual model is unconstrained, and the residual grows between them; the scheme-based residual couples consecutive instants and remains small throughout.}
    \label{fig:numerical_vs_exact}
\end{figure}

Figure~\ref{fig:numerical_vs_exact} summarizes performance using one-step residual--vs--time plots and phase-space trajectories from a single initial condition. 
The plots also illustrate the domain of influence for the loss terms.
To ensure that the intervals coupled by the scheme residual cover the full time domain, the grid spacing should satisfy $\Delta t \le h$. In this experiment, when $h=0.5$ and $N=11$ ($\Delta t=1>h$), the uncovered intervals lead to visible gaps in the learned dynamics. When $N\ge 21$, the union of intervals $[t,t+h]$ covers $[0,T]$, and the learned trajectories are correspondingly accurate.

The bottom row of Figure~\ref{fig:numerical_vs_exact} shows the corresponding results for the exact residual on the same fixed grids. Here the pointwise nature of the constraint is evident: even with $N=41$ fixed collocation points, the off-grid residuals reach $\mathcal{O}(1)$ near $t=0$, corrupting the early-time derivatives and, consequently, the phase-space trajectories; accurate results are obtained only when the time-collocation points are randomly resampled at each iteration (not shown). The scheme-based residual, by contrast, maintains $\mathcal{O}(10^{-2})$ off-grid residuals on the same fixed grids, consistent with the nonlocal coupling imposed by the scheme residual. We note that the two residual operators measure different quantities, so the residual magnitudes in the two rows are not directly comparable; the comparison rests on the resulting phase-space trajectories, obtained under identical sampling and training budgets.

\subsection{Critical points of the residual loss}\label{sec:critical-pts-of-R}

We next ask when a critical point of the residual functional necessarily has zero residual. The relevant condition for an implicit method is the local solvability of its update with respect to the new state. This condition should not be confused with numerical stability: a method may be unconditionally stable while the Jacobian of its implicit update becomes singular for particular step sizes or states.

Let $Q\in C^1(\mathbb{R}^{2d}\times\mathbb{R}^{2d};\mathbb{R}^{2d})$ define a one-step scheme by
\begin{equation}\label{eq:general-one-step-scheme}
    Q(u_{n+1},u_n)=0,
\end{equation}
and set
\begin{equation}\label{eq:general-one-step-residual}
    R_h[\Phi](u,t):=Q\lr{\Phi(u,t+h),\Phi(u,t)}.
\end{equation}
The result is not specific to Hamiltonian vector fields.

Here and below, a critical point means a critical point of the unrestricted function-space functional: the first variation vanishes for every admissible $C^1$ perturbation satisfying the initial condition.

\begin{thm}[Absence of spurious critical points for one-step schemes]\label{thm:general-one-step}
Let $\Omega\subset\mathbb R^{2d}$ be open, let $T\ge h$, let $\rho$ be a probability measure with full support on $\Omega$, and let $\gamma\in C([0,T])$ satisfy $\gamma(t)>0$ for every $t\in[0,T]$. Set
\[
    \mathcal R[\Phi]
    :=\frac12\int_\Omega\int_0^T
       \norm{R_h[\Phi](u,t)}^2\,\gamma(t)\,dt\,d\rho(u).
\]
Suppose $\Phi\in C^1(\Omega\times[0,T+h];\mathbb R^{2d})$ is a critical point among maps satisfying $\Phi(u,0)=u$. If the partial derivative of $Q$ with respect to its first argument,
\[
    \partial_1Q\lr{\Phi(u,t+h),\Phi(u,t)},
\]
is invertible for every $u\in\Omega$ and $t\in[0,T]$, then
\[
    R_h[\Phi](u,t)=0,
    \qquad u\in\Omega,\quad t\in[0,T].
\]
Hence every such critical point is a global minimizer. 
\end{thm}

\paragraph{Relation to neural parameterizations}
Theorem~\ref{thm:general-one-step} concerns the unrestricted function-space problem, whereas in computation we minimize $\mathcal R[\Phi_\theta]$ over a parameterized family $\mathcal M=\{\Phi_\theta:\theta\in\Theta\subset\mathbb R^P\}$ of neural networks. Let $\mathcal J(\theta):=\mathcal R[\Phi_\theta]$. By the chain rule,
\begin{equation}\label{eq:parameterized-first-variation}
    D\mathcal J(\theta)[\eta]
    =\delta\mathcal R[\Phi_\theta]\!\left[D_\theta\Phi_\theta[\eta]\right].
\end{equation}
Thus parameter-space stationarity requires the first variation to vanish only along directions in the network tangent space
\[
    T_{\Phi_\theta}\mathcal M:=\operatorname{Range}(D_\theta\Phi_\theta),
\]
rather than along all admissible function-space perturbations. The theorem therefore isolates the variational part of the problem: under its hypotheses, the residual functional itself has no nonzero-residual critical points. A nonzero-residual stationary point created after neural parameterization can occur only because the available tangent directions fail to detect a descent direction of the full functional. Consequently, when the network tangent space is sufficiently rich to approximate the relevant function-space descent directions, the theorem provides a natural guide to the parameterized problem, although it does not by itself establish a parameter-space landscape result.

For a finite collocation loss, the connection can be stated more concretely. Collect the network values entering the discrete loss into a vector $Z(\theta)$ and write the discrete objective as $\mathcal J(\theta)=\ell(Z(\theta))$. Then $\nabla_\theta \mathcal J=D_\theta Z(\theta)^T\nabla_Z \ell$. If $D_\theta Z(\theta)$ has full row rank, parameter-space stationarity implies stationarity with respect to all sampled output values. Thus, for the discretized problem, the remaining gap is explicitly a question of the rank and expressivity of the parameterization rather than of the residual functional alone.

The proof in~\ref{sec:proofs_critical_points} follows directly from the first variation. With $w(t)=R_h[\Phi](u,t)$, the optimality condition yields a terminal relation on $(T,T+h)$ and, after shifting the variation at $t+h$, a backward adjoint relation coupling $w(t-h)$ and $w(t)$. Invertibility of $\partial_1Q$ first forces the residual to vanish on the final interval and then propagates zero backward through the whole time domain.

\paragraph{Examples and schemes used here}
For an implicit--explicit representation
\[
Q(v,u)=\Phi_h^{\mathrm{Im}}(v)-\Phi_h^{\mathrm{Ex}}(u),
\]
Theorem~\ref{thm:general-one-step} requires $D\Phi_h^{\mathrm{Im}}$ to be invertible along the critical point. For an explicit one-step method, $\Phi_h^{\mathrm{Im}}=\mathrm{id}$, so the condition is automatic. For the Velocity Verlet scheme used in our residual experiments,
\[
\Phi_h^{\mathrm{Im}}(p,q)=
\begin{pmatrix}p-\tfrac h2F(q)\\ q\end{pmatrix},
\qquad
D\Phi_h^{\mathrm{Im}}(p,q)=
\begin{pmatrix}I&-\tfrac h2DF(q)\\0&I\end{pmatrix},
\]
which is invertible for every $h$ wherever $F$ is differentiable. For implicit midpoint,
\[
Q(v,u)=\frac{v-u}{h}-f\lr{\frac{u+v}{2}},
\qquad
\partial_1Q(v,u)=\frac1h\lr{I-\frac h2Df\lr{\frac{u+v}{2}}},
\]
so the relevant hypothesis is nonsingularity of the new-state equation. A small-step condition such as $\frac h2\|Df\|<1$ is sufficient but not necessary; this is a local solvability condition, not a CFL or stability restriction. The same distinction applies to other implicit schemes such as implicit Euler and the trapezoidal rule.

\begin{figure}[ht]
    \centering
    \begin{subfigure}[b]{0.42\linewidth}
        \centering
        \includegraphics[width=\linewidth,trim=90 10 190 10, clip]{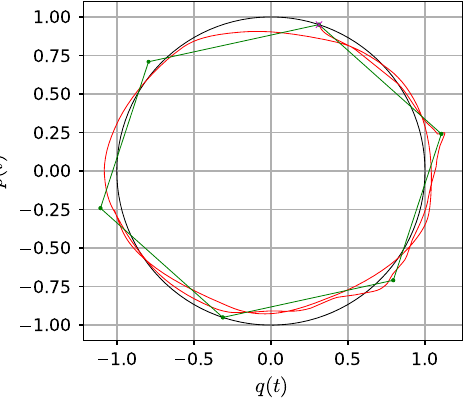}
        \caption{phase space}
    \end{subfigure}
    \begin{subfigure}[b]{0.42\linewidth}
        \centering
        \includegraphics[width=\linewidth,trim=0 10 180 10, clip]{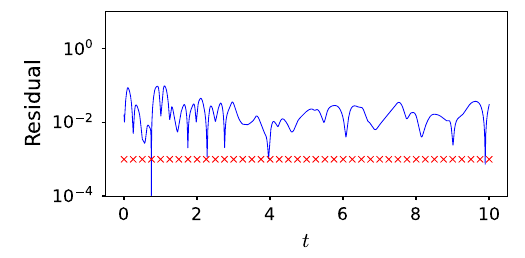}
        \caption{residual}
    \end{subfigure}
    \caption{(Harmonic oscillator) Effect of the coarse step $h=1$ used to define the residual, at $N=41$ collocation points. Black: exact flow. Green: the Velocity Verlet solution at $t=0,h,2h,\dots$. Red: the learned flow map. The learned map follows the numerical trajectory rather than the exact circle while the residual is still driven to a small value.}
    \label{fig:vv_stepsize}
\end{figure}

\paragraph{The learned map inherits the scheme, not the flow}
At zero residual, the learned map reproduces the selected numerical scheme rather than the exact flow. For a scheme written as an update map $u_{n+1}=\Psi_h(u_n)$ with residual
\[
R_h[\Phi](u,t)=\frac1h\bigl(\Phi(u,t+h)-\Psi_h(\Phi(u,t))\bigr),
\]
zero residual together with $\Phi(u,0)=u$ implies recursively that
\begin{equation}\label{eq:scheme-composition-target}
    \Phi(u,nh)=\Psi_h^n(u),\qquad n=0,1,2,\ldots.
\end{equation}
Thus the purpose of learning $\Phi$ is not to replace the underlying scheme by a larger numerical step, but to amortize the composition of many resolving steps into a map that can be evaluated directly at a coarse time increment. Theorem~\ref{thm:general-one-step} identifies this function-space target; it does not by itself guarantee its approximation or optimization within a neural parameterization. Figure~\ref{fig:vv_stepsize} makes this distinction visible by increasing the Velocity Verlet step to $h=1$: the learned trajectory follows the visibly distorted numerical trajectory rather than the exact circle. Thus reducing the residual cannot remove the discretization error of the underlying scheme.

\subsection{Error of the learned flow maps}

Having established conditions under which the critical points of the residual functional are global minimizers, we next quantify how a small scheme residual translates into flow-map error. For this estimate, write the selected one-step method as an update map $\Psi_h$ and use the normalized residual
\begin{equation}\label{eq:update-map-residual}
    R_h[\Phi](u,t)
    :=\frac1h\lr{\Phi(u,t+h)-\Psi_h(\Phi(u,t))}.
\end{equation}
This specialization is needed because a bound on a general implicit relation $Q(v,u)$ does not, by itself, provide a quantitative bound on the state-update error.

Assume (i) \emph{numerical accuracy}: on the relevant domain,
\[
    \norm{\Psi_h(v)-\phi_h(v)}\le C_{\rm loc}h^{p+1};
\]
(ii) \emph{training quality}: $\norm{R_h[\Phi](u,t)}\le\epsilon$ uniformly for $t\in[0,T]$, and the error on the initial interval satisfies
\[
    \delta_0:=\sup_{0\le s<h}\norm{\Phi(u,s)-\phi_s(u)}
    =\mathcal O(h^p)+\mathcal O(\epsilon);
\]
and (iii) the exact flow is Lipschitz with growth factor $e^{Lt}$. Then the following standard estimate holds.

\begin{prop}[Error of the learned flow map]\label{prop:flow-map-error}
Under assumptions (i)--(iii),
\begin{equation}
    \norm{\phi_t(u)-\Phi(u,t)}
    \le e^{Lt}\delta_0
    +\lr{\frac{C_{\rm loc}}{L}h^p+\frac{\epsilon}{L}}(e^{Lt}-1),
    \qquad t\in[0,T],
\end{equation}
with the usual continuous limit $(e^{Lt}-1)/L=t$ when $L=0$.
\end{prop}
\begin{proof}[Sketch]
Write $t=s+nh$ with $0\le s<h$ and start from the base error at time $s$. At each step, \eqref{eq:update-map-residual} contributes at most $h\epsilon$, while the numerical update differs from the exact time-$h$ flow by at most $C_{\rm loc}h^{p+1}$. Comparing exact flows from the learned and reference states contributes the factor $e^{Lh}$. Iterating
\[
    e_{k+1}\le e^{Lh}e_k+h\epsilon+C_{\rm loc}h^{p+1}
\]
and using $e^{Lh}-1\ge Lh$ gives the stated bound. This is the classical global-error argument for one-step methods (``Lady Windermere's fan''; see, e.g., \cite[Ch.~II]{hairer1993solving} and \cite{geometric}).
\end{proof}

The exponential factor in Proposition~\ref{prop:flow-map-error} is a worst-case estimate based on a global Lipschitz constant. For stiff Hamiltonian systems this can be excessively pessimistic because the large linear frequencies generate rotations rather than exponential growth. We now isolate this effect and obtain an amplification rate that is independent of the stiff frequencies.

\subsection{Stiffness-uniform error propagation for oscillatory Hamiltonians}\label{sec:ana_errorbound}
Consider
\begin{equation}\label{eq:osc-class}
    H(p,q)=\frac12\|p\|^2+\frac12q^T\Omega^2q+U(q),
    \qquad u=(p,q)\in\mathbb R^d\times\mathbb R^d,
\end{equation}
where $\Omega=\operatorname{diag}(\omega_1,\ldots,\omega_d)$ with $\omega_j\ge0$. The equations are $\dot p=-\Omega^2q-\nabla U(q)$ and $\dot q=p$.

\begin{assumption}[Uniform regularity on an energy band]\label{ass:system-class}
Fix $H_0$ and $\Delta H>0$, and let
\[
\mathcal E_{H_0}:=\{u:|H(u)-H_0|\le\Delta H\},
\qquad
\mathcal E'_{H_0}:=\{u:|H(u)-H_0|\le2\Delta H\}.
\]
There is a constant $c_U=c_U(H_0,\Delta H)$, independent of the frequencies $\omega_j$, such that on $\mathcal E'_{H_0}$
\[
    \|p\|+\|\Omega q\|+\|\nabla U(q)\|+\|\nabla^2U(q)\|\le c_U.
\]
\end{assumption}
This explicit energy-band bound is what the proof uses. It holds for the FPUT Hamiltonian because its nonlinear potential is nonnegative and coercive in the relevant combinations of coordinates; see~\ref{sec:hypotheses-verification}.

Let $\tilde\omega_j=\max\{\omega_j,1\}$ and $\tilde\Omega=\operatorname{diag}(\tilde\omega_1,\ldots,\tilde\omega_d)$. The weighted norm is
\begin{equation}\label{eq:W-norm-theory}
    \norm{(v_p,v_q)}_W:=\norm{(v_p,\tilde\Omega v_q)}_2,
\end{equation}
identical to the energy-balanced norm introduced in Section~\ref{sec:method_energy_norm}.

Let $u_{n+1}=\phi_{\Delta t}(u_n)$ be the exact sampled flow and $\tilde u_{n+1}=\Phi_{\Delta t}(\tilde u_n)$ the recursively applied learned map, with $\tilde u_0=u_0$.

\begin{assumption}[Local accuracy of the learned map]\label{ass:local-error}
On $\mathcal E'_{H_0}$,
\[
    \delta(u):=\Phi_{\Delta t}(u)-\phi_{\Delta t}(u),
    \qquad
    \sup_{u\in\mathcal E'_{H_0}}\norm{\delta(u)}_W\le\delta^*.
\]
\end{assumption}

\begin{thm}[Stiffness-uniform error propagation]\label{thm:linear-error-growth}
Let Assumptions~\ref{ass:system-class} and \ref{ass:local-error} hold and suppose $H(u_0)=H_0$. There exist constants $M>0$, $E_0>0$, and $C_H>0$, independent of $\Omega$, $n$, and $\delta^*$. Define
\[
    C:=M\Delta t,
    \qquad
    G_n(C):=\frac{e^{nC}-1}{e^C-1},
\]
with $G_n(0):=n$. For any $N$ satisfying
\[
    G_N(C)\,\delta^*\le E_0,
\]
the following estimates hold for every $n\le N$:
\begin{align}
    \norm{\tilde u_n-u_n}_W &\le G_n(C)\,\delta^*,\label{eq:stiffness-uniform-main}\\
    |\tilde u_{n,q_j}-u_{n,q_j}|&\le \tilde\omega_j^{-1}G_n(C)\,\delta^*,\label{eq:stiffness-uniform-position}\\
    |H(\tilde u_n)-H(u_n)|&\le C_H G_n(C)\,\delta^*.\label{eq:stiffness-uniform-energy}
\end{align}
In particular, the amplification rate $C/\Delta t=M$ is independent of the stiff frequencies. If also $NC\le c_0$ for a fixed constant $c_0$, then $G_n(C)\le e^{c_0}n$, and hence
\[
    \norm{\tilde u_n-u_n}_W\lesssim n\delta^*,
    \qquad
    |H(\tilde u_n)-H(u_n)|\lesssim n\delta^*.
\]
\end{thm}
Thus the theorem gives linear accumulation only while the slow-scale amplification remains bounded; it is not a global-in-time linear-growth result.

The proof is given in~\ref{sec:proof_global_error}. After rescaling the position variables by $\tilde\Omega$, the stiff linear part becomes a direct sum of rotations and therefore has unit operator norm. The remaining part of the linearized dynamics is bounded on the energy band by a constant independent of $\Omega$. This yields the one-step factor $e^C$ in \eqref{eq:stiffness-uniform-main}. The theorem therefore removes the spurious $\mathcal O(\max_j\omega_j)$ growth rate obtained from a direct Euclidean Lipschitz estimate; it does not rule out exponential sensitivity generated by the slow or nonlinear dynamics.

The result also explains the loss scaling. A Euclidean local-error bound $\|\delta(u)\|_2\le\delta_2^*$ only gives
\[
    \norm{\delta(u)}_W\le(\max_j\tilde\omega_j)\delta_2^*,
\]
which reintroduces the largest frequency into the global estimate. Controlling the local error directly in $\norm{\cdot}_W$ avoids this loss. For FPUT this is precisely the weighting used in \eqref{def:weighted-l2-norm}; the corresponding ablation in Section~\ref{sec:energy-loss-ablation} shows that the unweighted loss fails to resolve the fast energy exchange.

The FPUT Hamiltonian belongs to the class \eqref{eq:osc-class};~\ref{sec:hypotheses-verification} verifies Assumption~\ref{ass:system-class}. The NPCO problem has a slightly different form, but its linearization is already a direct sum of rotations and the same argument applies without the FPUT rescaling. The $\alpha$-particle problem is not covered because its stiff linear part depends on position through the magnetic field.

\section{Neural flow maps}\label{sec:method_arch}

We use two neural parameterizations. A fixed-timestep map $\Phi_{T_0}(u)$ approximates $\phi_{T_0}(u)$ for one prescribed increment $T_0$, while a variable-timestep map $\Phi(u,t)$ approximates $\phi_t(u)$ over a time interval. Both are designed for recursive application over horizons much longer than their training increment or training window.

\subsection{Fixed-timestep flow map $\Phi_{T_0}(u)$}

For the fixed-timestep map $\Phi_{T_0}:\mathbb{R}^{2d}\to\mathbb{R}^{2d}$, we use the state $u$ together with the vector field $f(u)$ as network inputs:
\begin{align}
    \Phi_{T_0}(u) := G(u, f(u)),
\end{align}
where $f: \mathbb{R}^{2d} \to \mathbb{R}^{2d}$ is the right-hand side of \eqref{eq:generic-Hamiltonian-system}, and $G: \mathbb{R}^{2d} \times \mathbb{R}^{2d} \to \mathbb{R}^{2d}$ is a neural network.  
We implement $G$ as a multilayer perceptron (MLP) with skip connections and gating operations following \cite{wang2024piratenets}.

\subsection{Variable-timestep flow maps $\Phi(u,t)$}\label{sec:Taylor-based-networks}

Because the learned map is applied recursively, its behavior near $t=0$ should be controlled explicitly. We therefore build the variable-timestep architecture from the Taylor expansion of the exact flow at $t=0$.
Specifically, for a chosen fixed consistency order $p\ge 0$ (a hyperparameter set before training), we design neural networks of the form 
$$\Phi(u,t) = u + \Psi^{(p)}(u,t)$$ 
whose time derivatives at $t=0$ match those of $\phi_t(u)$ up to order $p$.
Writing
\begin{align}
    \phi_t(u) = u + F(u,t), \quad 
    F(u,t) := \int_{0}^{t} f(\phi_\tau(u)) \, d\tau,
\end{align}
we require $\Psi^{(p)}(u,0) = 0$ and
\begin{align}\label{eq:consistency-derivative-matching}
    \left. \frac{\partial^k}{\partial t^k} \Psi^{(p)}(u,t) \right|_{t=0}
    = \left. \frac{\partial^k}{\partial t^k} F(u,t) \right|_{t=0}, \quad k=0,\dots,p.
\end{align}
Using the Lie derivative $D_f g(u) := g'(u) f(u)$, the Taylor expansion of $F$ in $t$ around $0$ reads
\begin{align}
    F(u,t) = \sum_{k=0}^\infty \frac{t^{k+1}}{(k+1)!} D_f^k f(u).
\end{align}
We approximate $F$ by truncating this series at order $p$ and adding a learnable remainder $\Delta^{(p+1)}_\theta$ that captures higher-order effects:
\begin{align}\label{eq:leading_order_taylor_generic}
    F(u,t) \approx \Psi^{(p)}(u,t):=\sum_{k=0}^{p-1} \frac{t^{k+1}}{(k+1)!} D_f^k f(u) + t^{p+1} \Delta^{(p+1)}_\theta(u, f(u), t).
\end{align}

The number of matched time derivatives at $t=0$ defines the consistency order of the neural flow map. Since $p$ is fixed, the Taylor terms are determined analytically from $f$; only $\Delta^{(p+1)}_\theta$ and the time-activation parameters (see below) are learned.

To curb the polynomial time-growth of the Taylor terms, we use simple but specialized activations: for example,
\begin{align}
    \Psi^{(0)}(u,t) &= \sigma(w_1 t) \Delta^{(1)}(u, f(u), t), \\
    \Psi^{(1)}(u,t) &= \frac{\sigma(w_1 t)}{w_1} f(u) + \frac{\sigma(w_1 t)}{w_1} \sigma(w_2 t) \Delta^{(2)}(u, f(u), t), \\
    \Psi^{(2)}(u,t) &= \frac{\sigma(w_1 t)}{w_1} f(u) + \frac12 \tfrac{\sigma(w_1 t)}{w_1} \tfrac{\sigma(w_2 t)}{w_2} f'(u) f(u) \nonumber \\
    &\quad + \frac{\sigma(w_1 t)}{w_1} \frac{\sigma(w_2 t)}{w_2} \sigma(w_3 t) \Delta^{(3)}(u, f(u), t), 
\end{align}
where $w_1, w_2, w_3>0$ are trainable parameters, and $\sigma$ satisfies
\[
\sigma(0)=0,\quad \sigma'(0)=1,\quad \sigma(\pm\infty) = \pm 1,
\]
(e.g., $\tanh$). In this way, $\sigma(w_i t)$ activates the Taylor terms over a learned time range: the expansion provides accuracy for small $t$, while $\Delta^{(p)}$ takes over for larger $t$. 

The following lemma confirms that the time activations do not destroy the built-in consistency.

\begin{lem}[Consistency of the activated architecture]\label{lem:consistency}
Let $f\in C^2$, let $\sigma\in C^3$ satisfy $\sigma(0)=0$ and $\sigma'(0)=1$, and let the learnable remainders $\Delta^{(k)}$ be $C^2$ in $t$. Then the activated flow maps $\Phi(u,t)=u+\Psi^{(p)}(u,t)$ defined above satisfy
\[
\left.\frac{\partial^k}{\partial t^k}\Phi(u,t)\right|_{t=0} = \left.\frac{\partial^k}{\partial t^k}\phi_t(u)\right|_{t=0}, \qquad k=0,\ldots,p,
\]
for $p=0$ and $p=1$; for $p=2$, the same holds provided additionally $\sigma''(0)=0$---a condition satisfied by every smooth odd activation, such as $\tanh$. The slow--fast variant \eqref{slow-fast-different-Taylor-orders} inherits the corresponding property componentwise.
\end{lem}

\begin{proof}
Direct differentiation at $t=0$. Every term containing a product of two (respectively three) activated factors $\sigma(w_it)$ vanishes at $t=0$ together with its first (respectively its first and second) derivatives, since each factor vanishes at $t=0$. For the leading terms, $\partial_t\lr{\sigma(w_1t)/w_1}f(u)\big|_{t=0}=f(u)$ matches $\partial_t\phi_t(u)|_{t=0}$; and $\partial_t^2\lr{\sigma(w_1t)/w_1}f(u)\big|_{t=0}=w_1\sigma''(0)f(u)$ vanishes when $\sigma''(0)=0$, leaving the second derivative of the second term, $f'(u)f(u)$, in agreement with $\partial_t^2\phi_t(u)|_{t=0}=D_ff(u)$.
\end{proof}

\paragraph{Varied orders for slow--fast systems}
In many physical systems (e.g., the FPUT chain in Section~\ref{subsec:FPU}), the state naturally decomposes into slow and fast components $u = (u_s, u_f)$, $f = (f_s, f_f)$. Using a single Taylor order $p$ for the entire state can be wasteful: the fast components oscillate rapidly and are better left to the learned remainder $\Delta$, while the slow components benefit from higher-order Taylor matching. We therefore define separate corrections $\Psi_s$ and $\Psi_f$ with possibly different Taylor orders $(p_s, p_f)$: 
\begin{equation}\label{slow-fast-different-Taylor-orders}
\begin{aligned}
    \Psi_s^{(p_s)}(u,t) &= \sum_{k=0}^{p_s-1} \frac{t^{k+1}}{(k+1)!} D_f^k f_s(u) + t^{p_s+1} \Delta_s(u, f(u), t), \\
    \Psi_f^{(p_f)}(u,t) &= \sum_{k=0}^{p_f-1} \frac{t^{k+1}}{(k+1)!} D_f^k f_f(u) + t^{p_f+1} \Delta_f(u, f(u), t),
\end{aligned}    
\end{equation}
and set $\Phi_s(u,t) = u_s + \Psi_s(u,t)$, $\Phi_f(u,t) = u_f + \Psi_f(u,t)$.  
Typically, $p_f=0$ for fast variables and $p_s>0$ for slow variables.  
Time activations can be applied here as in the generic case.
This decomposition requires coordinates in which the slow--fast splitting is explicit. When no such splitting is available, a uniform Taylor order can be used for all components.

\section{Algorithms for generating efficient training data}\label{sec:method_data}

\subsection{Constant energy manifolds and the microcanonical ensemble}
Since Hamiltonian flows preserve energy, it is natural to concentrate training points near dynamically relevant energy levels. For a target energy $H_0$, the microcanonical measure is formally
\begin{equation}
    \pi(dp\,dq)=Z^{-1}\,\delta\!\left(H(p,q)-H_0\right)\,dp\,dq,
    \label{eq:mcdensity}
\end{equation}
which, by the coarea formula, induces surface density proportional to $dS/\|\nabla H\|$ on $\Sigma_{H_0}=\{H=H_0\}$. This is the invariant measure obtained by restricting Liouville measure to the energy surface.

In practice we sample from several nearby energy levels rather than a single surface, because recursively applied learned maps may incur small energy errors. The experiments draw the energy from a narrow normal distribution centered at $H_0$ and then sample microcanonically on the corresponding level set. Compared with uniform sampling from a bounding box, this concentrates the sampling budget near the energies at which the learned map will be evaluated; Section~\ref{sec:res_phasespace} quantifies the resulting effect on training efficiency and long-time accuracy.

\subsection{Sampling algorithm: HMC-$H_0$}

Building on \cite{fang2024stabilization}, we use the following sampler for the microcanonical ensemble on $\Sigma_{H_0}$. Assume the Hamiltonian is separable, 
\begin{equation}
    H(p,q) = \frac{1}{2} p^T M^{-1} p + U(q),
\end{equation}
where $M$ is a diagonal mass matrix and $U(q)$ is a given potential energy. 

Algorithm~\ref{alg:hmc-H0} is closely related to Hamiltonian Monte Carlo (HMC) \cite{DUANE1987216,bou2018geometric}, but uses a different momentum refreshment. Standard HMC draws $p\sim\mathcal{N}(0,M)$ independently of the current position and targets a canonical (Gibbs) distribution. Here, conditional on the current position $q$, the momentum is sampled uniformly from the kinetic-energy surface
\[
    \frac12 p^TM^{-1}p=H_0-U(q),
\]
so that the refreshed state remains on the prescribed energy surface. The subsequent Hamiltonian integration time is randomized according to a distribution $\mu$; one possible choice is an exponential distribution with mean $\lambda$, as in \cite{bou2018geometric}.

HMC-$H_0$ is naturally parallel: independent chains can be initialized from different positions $q_0$.

The correctness of the exact-flow algorithm is summarized by the following proposition.

\begin{prop}[Invariance of the target measure under HMC-$H_0$]\label{prop:hmc-invariance}
Let $H(p,q)=\frac12 p^TM^{-1}p+U(q)$ and suppose the exact flow $\phi_t$ conserves $H$. Then the microcanonical measure
\begin{equation}\label{eq:hmc-target-measure}
    \pi(dp\,dq)=Z^{-1}\,\delta\!\lr{H(p,q)-H_0}\,dp\,dq
\end{equation}
is invariant under both steps of Algorithm~\ref{alg:hmc-H0}; hence $\pi$ is invariant for the algorithm.
\end{prop}

\begin{proof}
The refreshment step is a Gibbs update. Conditional on $q$, write $v=M^{-1/2}p$; the conditional density of \eqref{eq:hmc-target-measure} is proportional to
\[
    \delta\!\lr{\tfrac12\|v\|^2-(H_0-U(q))},
\]
so $v$ is uniform on the sphere of radius $\sqrt{2(H_0-U(q))}$, exactly as in Algorithm~\ref{alg:hmc-H0}. Resampling from this conditional leaves the joint measure invariant. The exact Hamiltonian flow preserves $dp\,dq$ by Liouville's theorem and preserves the level sets of $H$, so it also preserves $\pi$ for every fixed integration time. Averaging over an independently sampled $t\sim\mu$ therefore preserves $\pi$ as well.
\end{proof}

In practice, Algorithm~\ref{alg:hmc-H0} replaces the exact flow by an accurate symplectic numerical integrator. Without an accept--reject correction, the resulting Markov chain is therefore only an approximate microcanonical sampler: the exact invariance asserted in Proposition~\ref{prop:hmc-invariance} applies to the exact-flow version. In the experiments, the numerical integration error is chosen small relative to the width of the sampled energy band.

\begin{algorithm}
\caption{HMC-$H_0$}\label{alg:hmc-H0}
\begin{algorithmic}
\Require{initial position $q_0 \in \mathbb{R}^d$, mass matrix $M$, potential energy function $U$, target energy $H_0$ that satisfies $H_0-U(q_0)>0$, number of samples $N$, a distribution $\mu$ on $(0,\infty)$ with mean $\lambda$, reference flow map $\phi$}
\Ensure{phase space samples $\{ (p_{n}, q_{n}) \}_{n=1}^N$}
    \For{$n = 1:N$}
        \State Step 1: \textbf{Momentum refreshment.}
        \State sample ${\xi}_0 \sim \mathcal{U}(\mathbb{S}^{d-1})$
        \State $\xi = \sqrt{2(H_0-U(q_{n-1}))} M^{1/2} {\xi}_0$
        \State 
        \State Step 2: \textbf{Time integration.}
        \State sample $t \sim \mu$
        \State $p_{n}, q_{n} = \phi_{t}(\xi, q_{n-1})$
    \EndFor
\end{algorithmic}
\end{algorithm}

\section{Numerical experiments}\label{sec:dl-results}
 
We evaluate the proposed approach on three Hamiltonian systems: 
\begin{description}
    \item[Nearly-periodic coupled oscillators (NPCOs)] This is a weakly coupled slow-fast system studied in \cite{duruisseaux2023approximation}. We use this system to test various features of the proposed methodology.  
    \item[Fermi-Pasta-Ulam-Tsingou (FPUT) problem] This is a benchmark problem commonly used to compare geometric integrators for long-time simulations \cite{geometric}. With $\epsilon$ denoting the fast timescale, the standard long-time benchmark extends to horizons of order $\mathcal O(\epsilon^{-1})$. 
    \item[The $\alpha$-particle problem] This nonseparable, noncanonical system describes a charged particle in a magnetic-confinement model. It is derived and studied in \cite{burby2024non} in the context of stellarator modeling. We use this system to demonstrate that simulations with the proposed neural network, trained via residual minimization, exhibit long-time statistical patterns (Poincaré sections) that closely resemble the highly accurate reference solutions.
\end{description}

Section~\ref{subsec:runtime-comp} compares inference runtimes with conventional integrators.

\paragraph*{The HMC-$H_0$ datasets}
Some of the models reported below used datasets computed by Algorithm~\ref{alg:hmc-H0}. 
For each involved model, Algorithm~\ref{alg:hmc-H0} is applied to sample from the microcanonical ensemble on a finite number of energy level sets randomly selected according to a chosen normal distribution, centered at the target energy level, $H_0$, with the standard deviation set to $H_0/10.$ The mean time integration window in the algorithm is set to $\lambda=64$ for the NPCOs and $\lambda=0.5$ for the FPUT simulations. 
We refer to these samples as HMC-$H_0$ datasets.

\paragraph*{Error metrics}

In all the comparisons, we define the trajectory error between a predicted solution $\tilde{u}=(\tilde{p},\tilde{q})$ and the reference solution $u=(p,q)$ as  
\begin{align}\label{def:traj_err}
    \text{traj error} = \norm{\tilde{u}-u}_2,
\end{align}
and the energy error as  
\begin{align}\label{def:energy-error}
    \text{energy error} = \frac{\abs{H(\tilde{p},\tilde{q}) - H(p,q)} }{\abs{H(p,q)}}. 
\end{align}

Typically, the predicted solutions are computed by repeatedly applying a trained neural network or a reference numerical integrator. For example, when studying $\Phi(u_0,t),$
$\tilde u = \Phi^{(k)}(u_0, t),$ where
\begin{equation}\label{eq:recursion-of-a-flow-map}
 \Phi^{(k+1)}(u_0,t):=\Phi(\Phi^{(k)}(u_0,t),t),~~k=1,2,\cdots,~~~\Phi^{(1)}(u_0,t)=\Phi(u_0,t).
\end{equation}

In this section, readers will see additional superscripts and subscripts in the learned neural flow maps $\Phi$; they are used to specify the dependence on different training setups. The usage will be made clear in the relevant subsections.

We provide additional details on our numerical experiments in the Appendix. Ablation studies isolating the main design choices---the phase-space sampling distribution, the Taylor truncation order and the training horizon---are collected in Section~\ref{sec:ablations}, and the effect of the multi-step data loss in Section~\ref{sec:res_multistep}. 
For reproducibility, the Appendix records the model architectures, loss and sampling choices, training budgets, and the datasets used in the sampling ablations (Tables~\ref{tab:training_data_ablation}--\ref{tab:flow_map_training}).

\subsection{The nearly periodic coupled oscillators (NPCO) problem}\label{subsec:NCPO}

The system is defined by the Hamiltonian 
\begin{align}
    H_\epsilon (p_1,p_2,q_1,q_2) = \frac{1}{2} (q_1^2 + p_1^2 ) + \frac{1}{2} \epsilon (q_2^2 + p_2^2 )  + \epsilon U(q_1,q_2),
    \label{eq:H_npco}
\end{align}
where the coupling potential is defined by 
\begin{align}
    U(q_1,q_2) = q_1 q_2 \sin(2q_1 + 2q_2).
\end{align}
Here, the parameter $\epsilon$ controls both the coupling strength and separation of timescales between the oscillators. In the regime $\epsilon\ll 1$, the fast oscillator, with variables $(p_1,q_1)$, undergoes simple, nearly circular oscillations. In contrast, the slow oscillator, $(p_2,q_2)$, evolves on a timescale that is $\epsilon^{-1}$ times longer and displays intricate dynamical behavior in the phase space. 

We train a variable-timestep flow map $\Phi(u,t)$ over the time interval $[0,5]$ using collocation points drawn from a bounding box (Box) phase space distribution. As seen in Table~\ref{tab:flow_map_training}, the final train/test residual losses are both $1.5\times10^{-7}$, which suggests excellent agreement between the learned flow map and the underlying velocity Verlet scheme. 

\begin{figure}
    \centering
    \begin{subfigure}[b]{0.24\linewidth}
      \includegraphics[width=\linewidth,trim=170 0 0 15, clip]
      {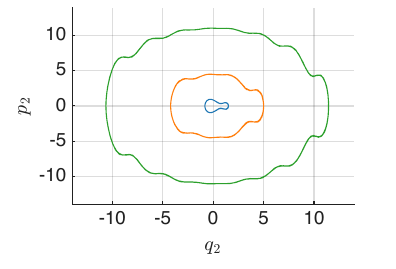}
        \caption{$\Phi^{(k)}(u_0,5)$}
    \end{subfigure}
    \begin{subfigure}[b]{0.24\linewidth}
      \includegraphics[width=\linewidth,trim=170 0 0 15, clip]
      {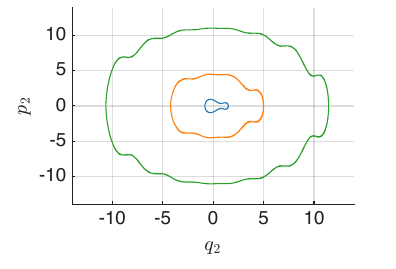}
        \caption{$\phi_{k\Delta t}(u_0)$}
    \end{subfigure}
    \begin{subfigure}[b]{0.24\linewidth}
      \includegraphics[width=\linewidth,trim=170 0 0 15, clip]
      {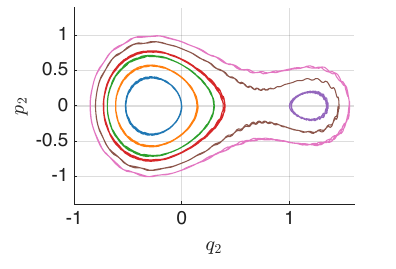}
        \caption{$\Phi^{(k)}(u_0,5)$}
    \end{subfigure}  
    \begin{subfigure}[b]{0.24\linewidth}        \includegraphics[width=\linewidth,trim=170 0 0 15, clip]
    {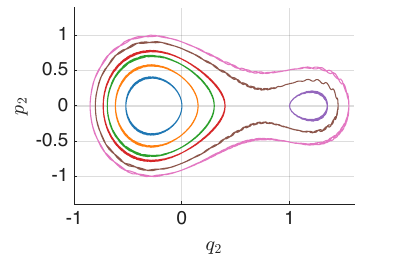}
        \caption{$\phi_{k\Delta t}(u_0)$}
    \end{subfigure}
    \caption{(NPCOs problem): Trajectories of the slow oscillator for $t\in [0,2000]$ generated by repeatedly applying the neural network-based flow map $\Phi(\cdot,\Delta t)$, $\Delta t=5$, and the reference flow map $\phi_t$. 
    (a),(b): three trajectories of the slow oscillator with the same energy $H_0=1.13$.   
    (c),(d): trajectories on different energy level sets in the inner region (where the blue curves in the top row occupy). }
    \label{fig:npco_nn_numerical}
\end{figure}

Figure~\ref{fig:npco_nn_numerical} compares long-time slow oscillator trajectories (up to $T=2000$) computed by $\Phi^{(k)}(u_0,\Delta t)$, in increments of $\Delta t=5$, to those of the reference flow $\phi_{k\Delta t}(u_0)$. The agreement persists after hundreds of successive applications of the learned map, demonstrating long-term stability and accuracy.


\subsection{The Fermi-Pasta-Ulam-Tsingou Problem}\label{subsec:FPU}

The Fermi–Pasta–Ulam–Tsingou (FPUT) problem describes a chain of masses connected by stiff linear and weakly nonlinear springs \cite{fermi1955studies}.
We consider $2m$ particles with displacements $q_i$ ($q_0=q_{2m+1}=0$) and momenta $p_i=dq_i/dt$. The system is described by the Hamiltonian: 
\begin{align*}
H(p,q)=&\frac{1}{2}\sum_{i=1}^m (p_{2i-1}^2+p_{2i}^2)
+ \frac{\omega^2}{4}\sum_{i=1}^m (q_{2i}-q_{2i-1})^2\\
&+ \sum_{i=0}^m (q_{2i+1}-q_{2i})^4,\qquad\omega\gg 1.
\end{align*}

Introducing slow variables $(x_{s,i},y_{s,i})$ and fast variables $(x_{f,i},y_{f,i})$ via 
    \begin{align*}
        x_{s,i} := \frac{1}{\sqrt{2}} \lr{q_{2i} + q_{2i-1}}, \quad y_{s,i} := \frac{1}{\sqrt{2}} \lr{p_{2i} + p_{2i-1}},
    \end{align*}
    \begin{align*}
        x_{f,i} := \frac{1}{\sqrt{2}} \lr{q_{2i} - q_{2i-1}}, \quad y_{f,i} := \frac{1}{\sqrt{2}} \lr{p_{2i} - p_{2i-1}},
    \end{align*}
for $i=1,\ldots,m$, 
the Hamiltonian separates into slow and oscillatory parts:
\begin{align}\label{eq:FPUT-H(y,x)}
    H(y,x) = & \frac{1}{4}\lr{x_{s,1}-x_{f,1}}^4+\frac{1}{4}\lr{x_{s,m}+x_{f,m}}^4\\
    +\frac{1}{2} \sum_{i=1}^{m} &(y_{s,i}^2 + y_{f,i}^2) + \frac{\omega^2}{2} \sum_{i=1}^{m} x_{f,i}^2 
    + \frac{1}{4}\sum_{i=1}^{m-1} \lr{x_{s,i+1}-x_{f,i+1}-x_{s,i}-x_{f,i}}^4. \nonumber
\end{align}

For this system, geometric averaging theory \cite[Ch.~XIII]{geometric} implies that the sum of the stiff-spring energies $I_j= \tfrac{1}{2}(y_{f,j}^2+\omega^2 x_{f,j}^2)$ is nearly conserved; i.e.,
\begin{align}
I=\sum_{j=1}^m I_j=I(y(0),x(0))+\mathcal{O}(\omega^{-1}).
\end{align}

The computational challenge is to resolve the nontrivial energy transfer among the stiff springs over an $\mathcal O(\omega)$ horizon while maintaining their nearly conserved total energy. Because the fast period is $\mathcal O(\omega^{-1})$, a direct method that resolves every oscillation requires at least $\mathcal O(\omega^2)$ steps over such a horizon.

\subsubsection{Effect of the energy-balanced loss}\label{sec:energy-loss-ablation}

For the FPUT variables $u=(y_s,x_s,y_f,x_f)$, the general weighting \eqref{def:weighted-l2-norm} reduces to
\[
    W=\operatorname{diag}(I_{3m},\omega I_m).
\]
Thus errors in the fast positions are multiplied by $\omega$ in both the data and residual losses. This is essential when $\omega\gg1$: on an $\mathcal O(1)$ energy surface the fast displacements are only $\mathcal O(\omega^{-1})$, so an unweighted mean-squared loss can be small while the corresponding stiff-spring energy is inaccurate. Theorem~\ref{thm:linear-error-growth} gives the complementary propagation result: local errors controlled in this weighted norm have an amplification rate independent of $\omega$.

The following comparison isolates the effect of the weighting; all other training choices are held fixed. The remaining FPUT experiments use the energy-balanced losses unless stated otherwise.

\paragraph{Quick comparisons}

We first compare two fixed-timestep flow maps trained under identical conditions 
(same dataset, architecture, and optimization setup) but differing in their loss metrics:
\begin{itemize}
    \item $\Phi_{T_0}^{\mathrm{MSE}}(u)$: trained by minimizing the standard mean-squared error between predicted and reference states;
    \item $\Phi_{T_0}^{\mathrm{EB}}(u)$: trained using the energy-balanced loss defined in \eqref{def:weighted-l2-norm}--\eqref{def:energy-balanced-residual-loss}.
\end{itemize}

\begin{figure}
    \centering
    \begin{subfigure}{0.45\linewidth}
        \includegraphics[valign=c, width=\linewidth, trim=0 30 0 0, clip]{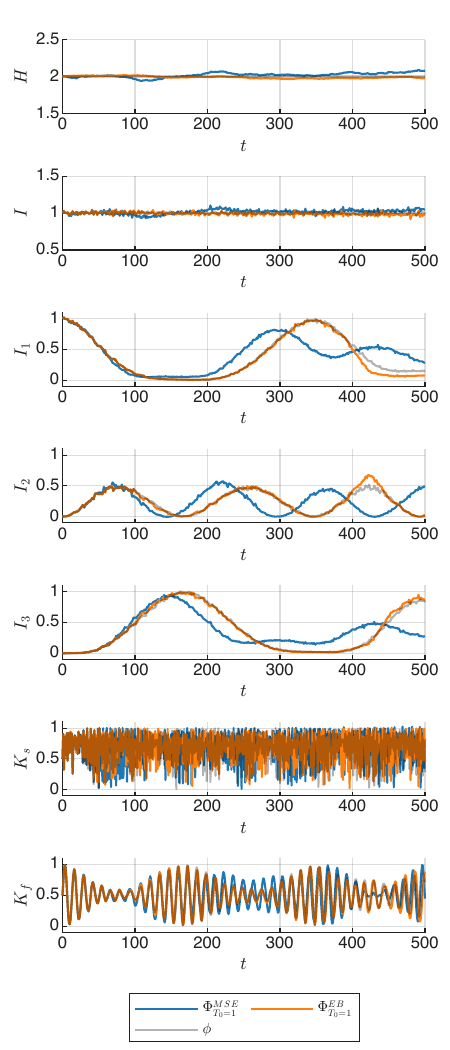}
    \end{subfigure}
    \begin{subfigure}{0.45\linewidth}
        \includegraphics[valign=c, width=\linewidth,trim=0 0 0 80, clip]  {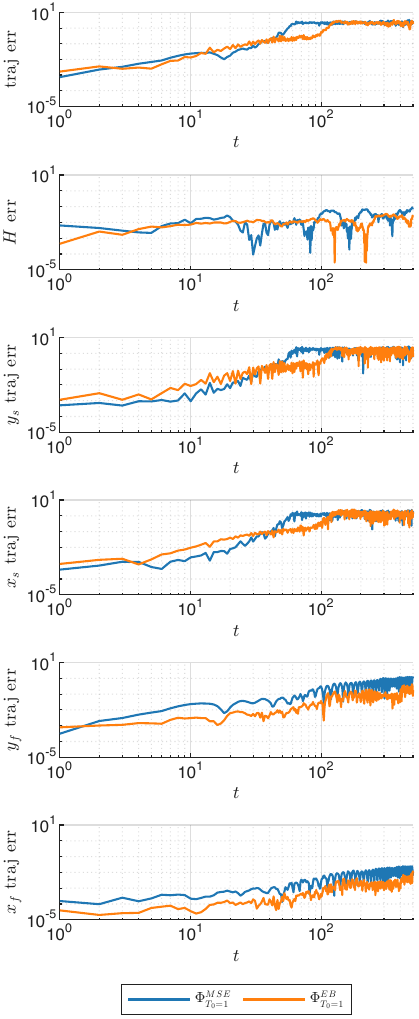}  
  
    \end{subfigure}
    \caption{(FPUT problem $\omega=50$) Errors and energy profiles of long-time trajectories over the interval $[0,500]$ generated by the fixed-timestep flow maps $\Phi_{T_0=1}$ trained with $\ell_2$ data loss and energy-balanced data loss.}
    \label{fig:fput_energybasednorm_dataloss}
\end{figure}

Figure~\ref{fig:fput_energybasednorm_dataloss} compares their long-time trajectories over $[0,500]$. 
The EB model $\Phi_{T_0}^{\mathrm{EB}}$ produces consistently more accurate fast components and, 
despite slightly larger initial errors in the slow components ($t\in[0,50]$), exhibits much slower error growth thereafter. 
The EB formulation also yields smaller overall energy errors, leading to improved conservation of total energy, 
the stiff-spring energy $I$, and the corresponding energy exchange dynamics.

Next, we examine two variable-timestep flow maps trained via a numerical residual loss:
\begin{itemize}
    \item $\Phi^{\mathrm{MSE}}(u,t)$: trained by minimizing the MSE of the residual;
    \item $\Phi^{\mathrm{EB}}(u,t)$: trained by minimizing the residual measured in the EB norm.
\end{itemize}
Time collocation points are sampled uniformly from $[0,0.125]$, and phase-space points from the HMC-$H_0$ dataset (with $\lambda=0.5)$. 
Long-time trajectories are obtained by iterating the learned maps 
$(\Phi^{\mathrm{MSE}})^{(k)}(u,t=0.0625)$ and $(\Phi^{\mathrm{EB}})^{(k)}(u,t=0.0625)$ for $kt\in[0,100]$. 

\begin{figure}
    \centering
    \begin{subfigure}{0.45\linewidth}
        \includegraphics[valign=c, width=\linewidth, trim=0 30 0 0, clip]
        {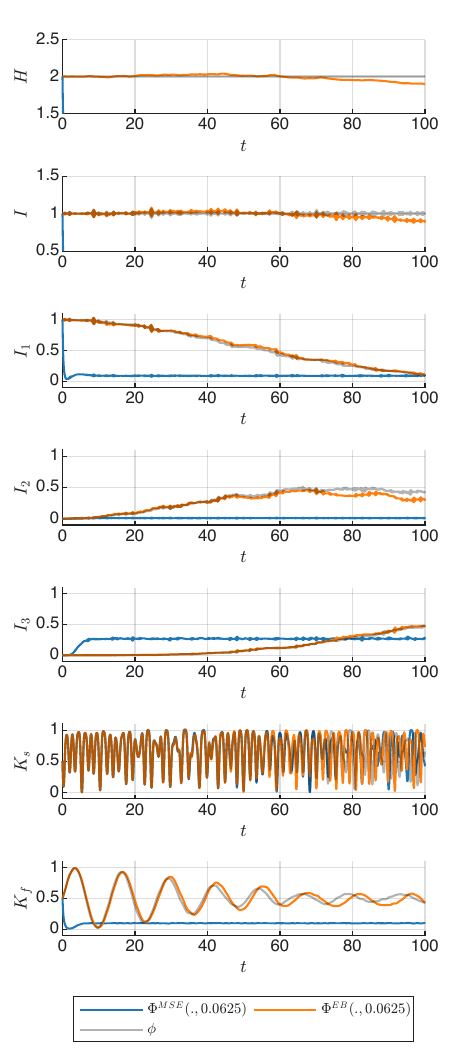}
    \end{subfigure}
    \begin{subfigure}{0.45\linewidth}
        \includegraphics[valign=c, width=\linewidth,trim=0 0 0 80, clip]{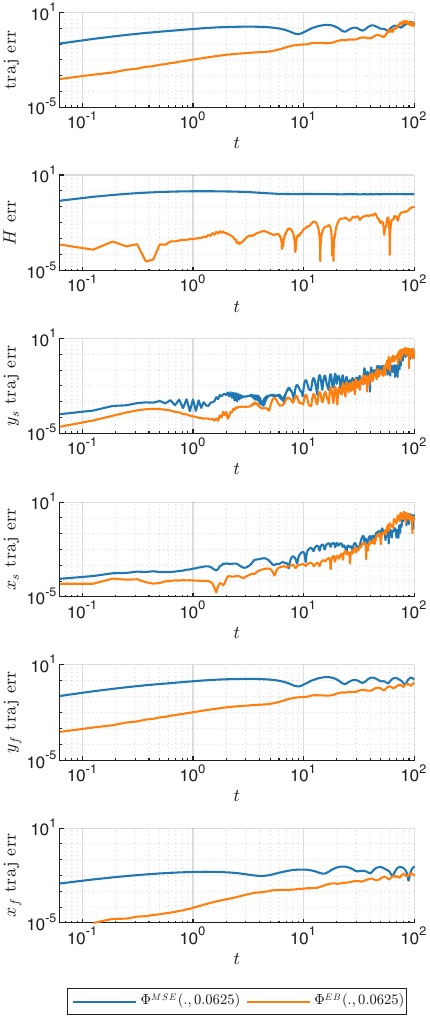}
    \end{subfigure}
    \caption{(FPUT problem $\omega=50$) Errors and energy profiles of long-time trajectories over the interval $[0,100]$ generated by recursively applying the flow maps with increment $t=0.0625$. $\Phi(u,t)$ is trained with $\ell_2$ residual loss and energy-balanced residual loss.}
    \label{fig:fput_energybasednorm_residualloss}
\end{figure}

Figure~\ref{fig:fput_energybasednorm_residualloss} shows that the MSE model fails to capture the fast oscillations: 
the fast kinetic energy $K_f$ collapses to zero almost immediately, indicating that the oscillatory modes are not learned.
In contrast, the EB model accurately preserves both the fast and slow dynamics, 
maintains the fast kinetic energy over the entire time interval, and achieves superior long-time accuracy.

\subsubsection{Further numerical results}
\label{sec:res_modelproblems}

We investigate two parameter regimes of the FPUT problem, $\omega=50$ and $\omega=300$. The FPUT system exhibits distinct dynamical features on different timescales: $\mathcal{O}(\omega^{-1})$, $\mathcal{O}(\omega^{0})$, $\mathcal{O}(\omega)$. Our primary goal here is to learn a valid flow map with a step size of $\mathcal{O}(\omega^{0})$ that is capable of predicting accurate long-time behavior on the scale $\mathcal{O}(\omega)$. The regime $\omega=300$ is especially challenging due to the significant separation of scales. We focus on reproducing conserved or nearly conserved quantities for long-time integration, i.e., the total energy $H$ and the total energy $I$ of the stiff springs.

\begin{figure}
    \centering
    \begin{subfigure}[b]{0.32\linewidth}
        \includegraphics[width=\linewidth]
        {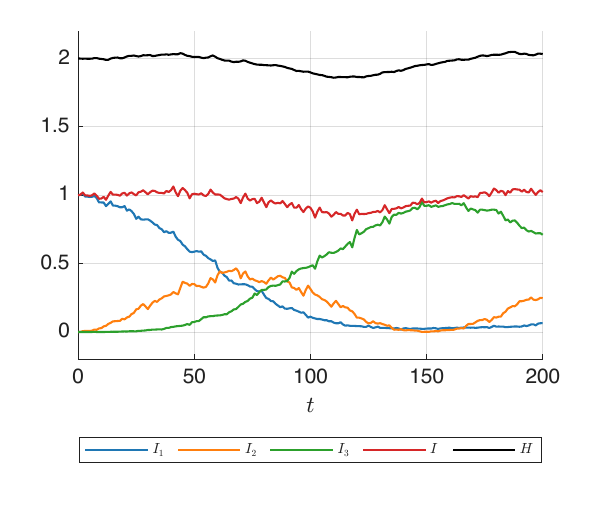}
        \caption{$\Phi(.,0.125)$}
    \end{subfigure}
    \begin{subfigure}[b]{0.32\linewidth}
        \includegraphics[width=\linewidth]
        {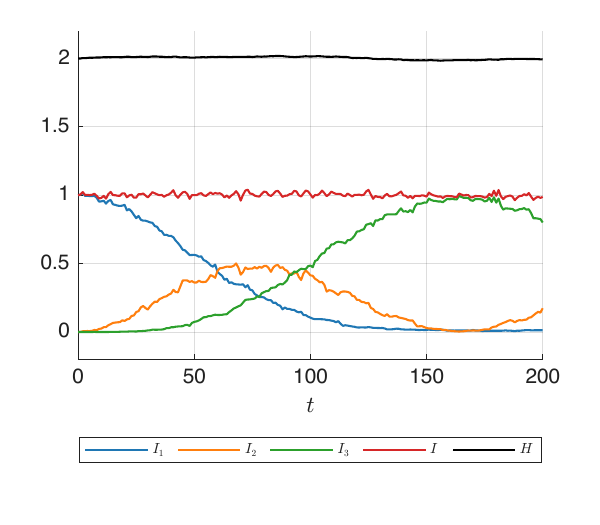}
        \caption{$\Phi_{T_0=1}$}
    \end{subfigure}
    \begin{subfigure}[b]{0.32\linewidth}
        \includegraphics[width=\linewidth]
        {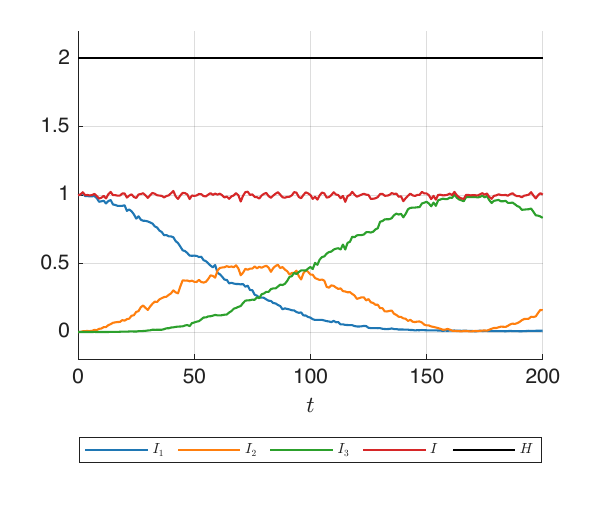}
        \caption{$\phi$}
    \end{subfigure}
    \caption{(FPUT problem $\omega=50$) Long-time energy profiles generated by neural network-based flow maps $\Phi(.,0.125)$, $\Phi_{T_0=1}$ and the reference flow map $\phi$.}
    \label{fig:fput_nn_numerical_omega_50}
\end{figure}

\paragraph{Case $\omega=50$} We train two distinct models: a variable-timestep flow map $\Phi$ over $[0,0.125]$ and a fixed-timestep flow map $\Phi_{T_0=1}$. Figure~\ref{fig:fput_nn_numerical_omega_50} shows the long-time energy profiles of the stiff springs. The reference solution $\phi$ exhibits characteristic exchange of energy among the stiff springs. The fixed-timestep map $\Phi_{T_0=1}$ closely tracks the reference energy exchange. The variable-timestep map $\Phi(\cdot,0.125)$ captures the same qualitative behavior but accumulates a somewhat larger phase error. Both remain stable over the reported horizon.

\begin{figure}
    \centering
    \begin{subfigure}[b]{0.32\linewidth}
        \includegraphics[width=\linewidth]{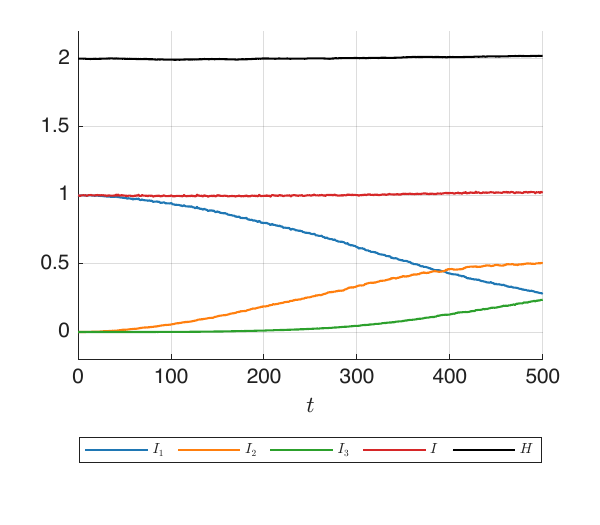}
        \caption{$\Phi_{T_0=1}$}
    \end{subfigure}
    \begin{subfigure}[b]{0.32\linewidth}
        \includegraphics[width=\linewidth]
        {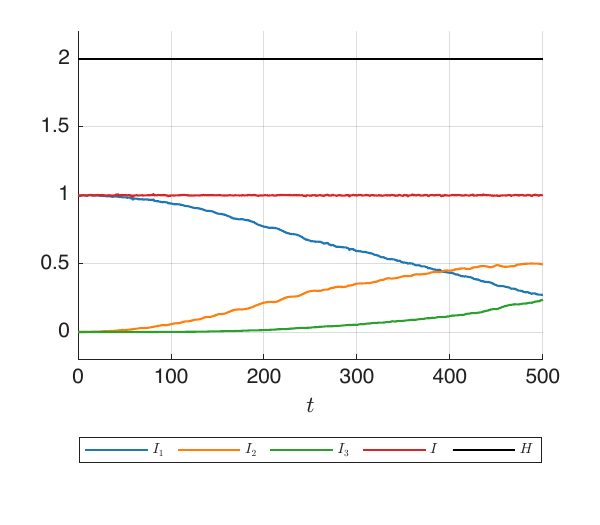}
        \caption{$\phi$}
    \end{subfigure}
    \caption{(FPUT problem $\omega=300$) Long-time energy profiles generated by the neural network-based flow map $\Phi_{T_0=1}$ and the reference flow map $\phi$.}
    \label{fig:fput_nn_numerical_omega_300}
\end{figure}

\paragraph{Case $\omega=300$} For this regime, we train only a fixed-timestep flow map $\Phi_{T_0=1}$, matching data evolved for multiple steps in \eqref{def:data-loss} with $S=5$. The trained model attains small train and test losses (Table~\ref{tab:flow_map_training}). Figure~\ref{fig:fput_nn_numerical_omega_300} compares its energy profiles with the reference solution; on the plotted scale, no systematic drift is visible in the total energy or in the distribution of energy among the stiff springs. As discussed later (Section~\ref{sec:runtime}), for large $\omega$ values, the neural network-based flow map can offer a runtime advantage over explicit solvers that must take very small time steps to resolve the fast scales.

\subsection{The $\alpha$-particle problem}\label{sec:runtime}

The dynamics of fusion-born $\alpha$-particles play a key role in magnetic confinement devices such as stellarators. We adopt the reduced model from \cite{burby2024non}, where the motion in a magnetic field $\bm{B}(x,y,z)=B(x,y)\bm{e}_z$ is described by
\begin{align}
\dot{v}_x = B(x,y)v_y, \quad \dot{v}_y=-B(x,y)v_x, \quad 
\dot{x} = \epsilon v_x, \quad \dot{y}=\epsilon v_y,
\end{align}
with $(v_x,v_y,x,y)\in\mathbb{R}^4$. This arises from a symmetry reduction of the six-dimensional Lorentz system under translation in $z$. The scale parameter $\epsilon=\rho/L$ (gyroradius to magnetic length scale) controls the degree of time-scale separation.

Expressing $v_x=r\cos\theta$, $v_y=r\sin\theta$, one finds $\dot{r}=0$ and $\dot{\theta}=-B(x,y)$, so $r$ is conserved while the gyromotion has frequency $B(x,y)$. The associated timescales are $T_{\text{gyro}}\sim 1/B(x,y)$ for the fast rotation and $T_{\text{drift}}\sim 1/(\epsilon B(x,y))$ for the slow guiding-center drift. Thus, the system naturally admits a multiscale decomposition, with ${\theta=0}$ providing a convenient Poincaré section for capturing the drift dynamics.
For magnetic fields of the form
\begin{align}
B(x,y)=B_0+a_1\cos(k_1 x+k_2 y)+a_2\cos(k_{3}x+k_{4}y),
\end{align}
the dynamics become increasingly chaotic as $\epsilon$ grows. Classical guiding center theory is only valid for $\epsilon\lesssim0.15$, so accurate long-time prediction of fusion-born $\alpha$-particles requires resolving the full reduced dynamics. However, the step size of an explicit method is restricted by the fast gyromotion, and the nonseparable, noncanonical Hamiltonian structure excludes standard symplectic schemes (e.g., the Verlet scheme). This setting motivates the design of structure-preserving integrators and learned flow-map surrogates that can uniformly reproduce the system’s multiscale features across $\epsilon$ and magnetic field configurations. 

To capture the parameter dependence on $\epsilon$, which controls the transition from regular to chaotic orbits, we include $\epsilon$ as an additional network input.
Writing $u=(v_x,v_y,x,y)$ and $\Phi(u,t) = (\Phi_v(u,t), \Phi_x(u,t))$, we directly enforce this property by letting 
\begin{align*}
    \Phi(u,t, \epsilon) &:= \lr{\frac{\| v\|}{\| \Phi_v(u,t, \epsilon)\|}\Phi_v(u,t, \epsilon), \Phi_x(u,t, \epsilon)}.
\end{align*}

We then train a variable-timestep flow map $\Phi$ on the interval $[0,5]$ via a numerical residual loss. For the phase space collocation points, we sample     $(x,y)\sim \mathcal{U}(\Omega_{xy})$ where $\Omega_{xy}\in{\mathbb{R}^2}$ covers several periodic cells of the magnetic field, and $(v_x,v_y)$ are sampled randomly from a thin shell with radius $r\in[\sqrt{2}-0.3, \sqrt{2}+0.3]$. The $\epsilon$ collocation points are sampled according to $\epsilon \sim \mathcal{U}([0.05,0.4])$. 

\begin{figure}
    \centering
    \begin{subfigure}[b]{0.35\linewidth}
        \includegraphics[width=\linewidth]
        {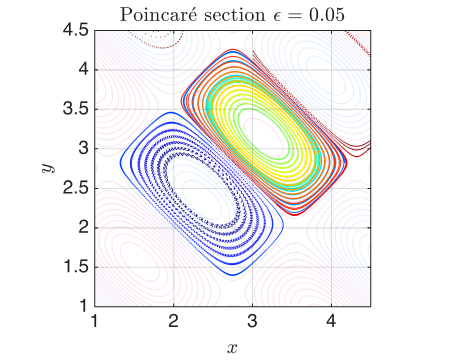}
        \caption{$\Phi(.,t)\;t\in[0,2]$}
    \end{subfigure}
    \begin{subfigure}[b]{0.35\linewidth}
        \includegraphics[width=\linewidth]
        {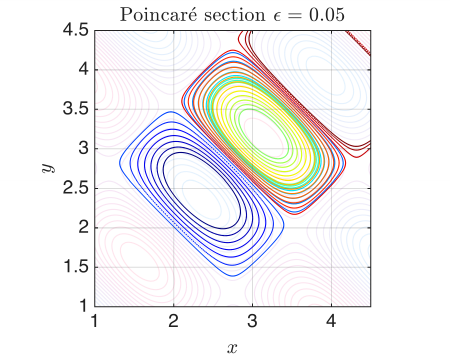}
        \caption{$\phi$}
    \end{subfigure}
    \par\bigskip
    \begin{subfigure}[b]{0.35\linewidth}
        \includegraphics[width=\linewidth]
        {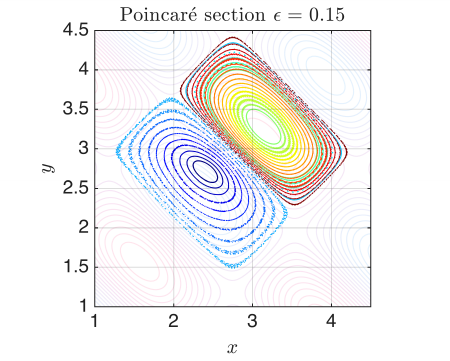}
        \caption{$\Phi(.,t)\;t\in[0,2]$}
    \end{subfigure}
    \begin{subfigure}[b]{0.35\linewidth}
        \includegraphics[width=\linewidth]
        {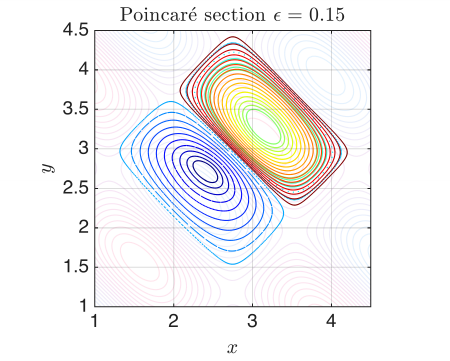}
        \caption{$\phi$}
    \end{subfigure}
    \par\bigskip
    \begin{subfigure}[b]{0.35\linewidth}
        \includegraphics[width=\linewidth]
        {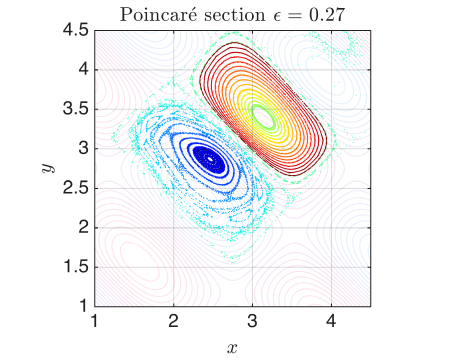}
        \caption{$\Phi(.,t)\;t\in[0,2]$}
    \end{subfigure}
    \begin{subfigure}[b]{0.35\linewidth}
        \includegraphics[width=\linewidth]
        {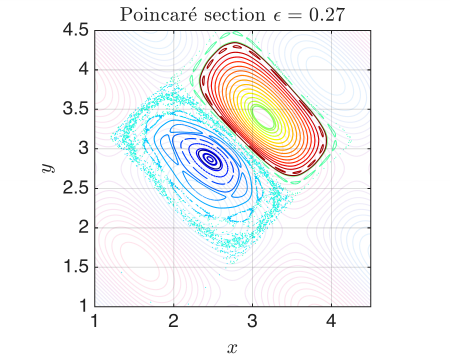}
        \caption{$\phi$}
    \end{subfigure}
    \caption{(The $\alpha$-particle problem) Poincaré sections at several values of $\epsilon$ generated by the neural network-based flow map $\Phi$ and the reference flow map $\phi$.}
    \label{fig:alphaparticle_nn_numerical}
\end{figure}

Figure~\ref{fig:alphaparticle_nn_numerical} shows Poincaré sections for several values of $\epsilon$, generated by the learned flow map $\Phi$ and the reference solution $\phi$. Over a long-time horizon $T=5000$, the learned solutions reproduce the same structures in the Poincaré sections, including the transition from more regular orbits at smaller $\epsilon$ to more chaotic trajectories at larger $\epsilon$. These results highlight the flow map's ability to capture the dependence on $\epsilon$ and to remain accurate across thousands of repeated applications.

\subsection{Ablation of the principal design choices}\label{sec:ablations}
The framework combines three choices that are logically independent of one another: where in phase space the training points are drawn, what is fitted and in which norm, and which function class carries the fit. The energy-balanced norm was isolated in Section~\ref{sec:energy-loss-ablation}, where removing it caused the fast kinetic energy to collapse. This section isolates the remaining two, together with the length of the training horizon and the number of steps matched by the data loss, all of which are free parameters of every experiment reported above.

\subsubsection{The training distribution}\label{sec:res_phasespace}
For the NPCO problem we train the fixed-timestep map $\Phi_{T_0=5}$ with a one-step data loss on two datasets that differ only in how the states are drawn: the HMC-$H_0$ dataset of Algorithm~\ref{alg:hmc-H0}, concentrated near $H_0=1.13$, and a dataset drawn uniformly from a bounding box in $\mathbb R^4$. The dataset-generation parameters are recorded in Table~\ref{tab:training_data_ablation} in the Appendix.

Figure~\ref{fig:npco_data}~(a) plots the test loss against the number of training samples. The model trained on HMC-$H_0$ data reaches a given test loss with substantially fewer samples: concentrating the samples on the surface where the map will actually be evaluated buys accuracy per sample that a wider, less relevant dataset does not. At a fixed budget of $24{,}000$ samples, the long-time trajectories over $[0,2000]$ in Figure~\ref{fig:npco_data}~(b)--(d) show the same ordering qualitatively---the HMC-$H_0$ model stays concentrated and close to the reference, while the bounding-box model drifts off the reference torus.

\begin{figure}[ht]
    \centering
    \begin{subfigure}[b]{0.5\linewidth}
        \includegraphics[width=\linewidth]{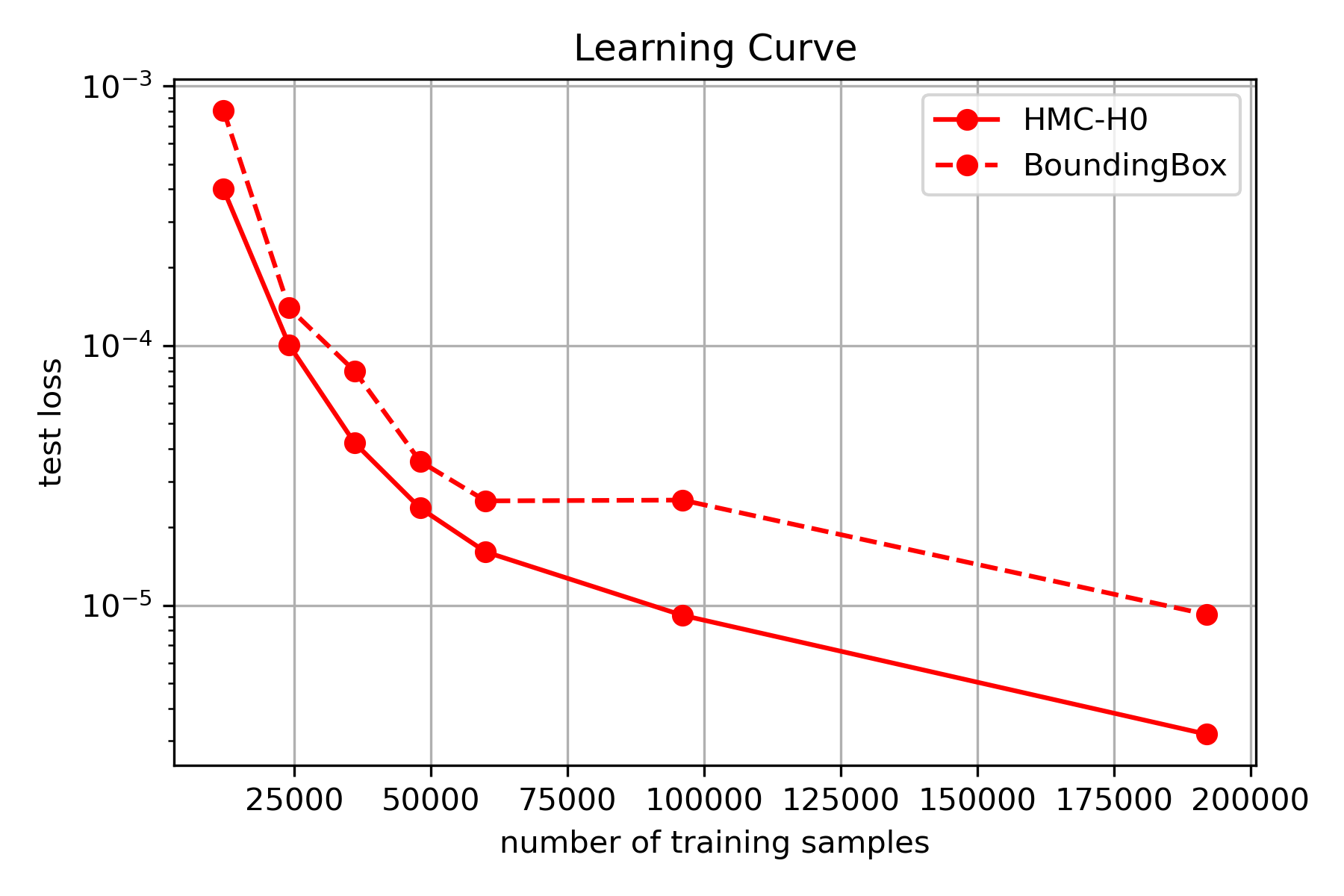}
        \caption{test loss vs.\ number of samples}
    \end{subfigure}
    \par\bigskip
    \begin{subfigure}[b]{0.32\linewidth}
        \includegraphics[width=\linewidth]{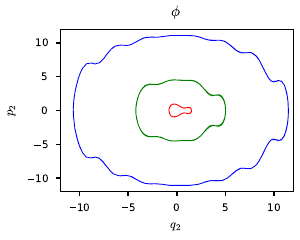}
        \caption{reference}
    \end{subfigure}
    \begin{subfigure}[b]{0.32\linewidth}
        \includegraphics[width=\linewidth]{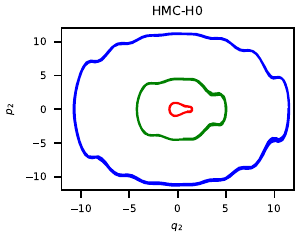}
        \caption{HMC-$H_0$ data}
    \end{subfigure}
    \begin{subfigure}[b]{0.32\linewidth}
        \includegraphics[width=\linewidth]{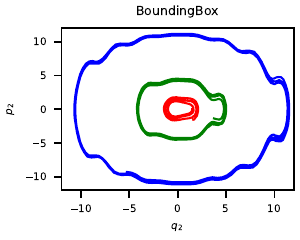}
        \caption{bounding-box data}
    \end{subfigure}
    \caption{(NPCO problem) Effect of the training distribution. (a) Test loss of $\Phi_{T_0=5}$ against the number of training samples. (b)--(d) Slow-oscillator phase-space plots of long-time trajectories over $[0,2000]$ generated by $\Phi_{T_0=5}$ trained on $24{,}000$ HMC-$H_0$ samples and on $24{,}000$ bounding-box samples.}
    \label{fig:npco_data}
\end{figure}

\begin{figure}[ht]
    \centering
    \begin{subfigure}[b]{0.45\linewidth}
        \includegraphics[width=\linewidth]{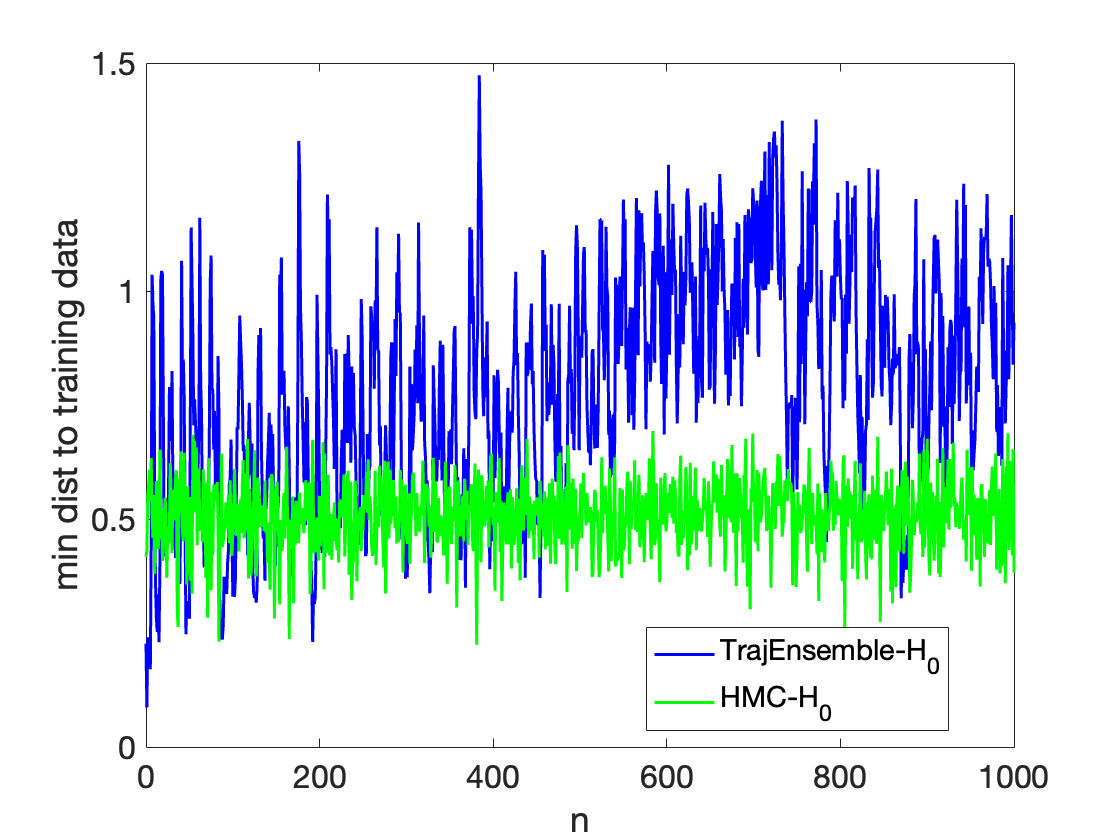}
        \caption{distance to the training set}
    \end{subfigure}
    \begin{subfigure}[b]{0.45\linewidth}
        \includegraphics[width=\linewidth]{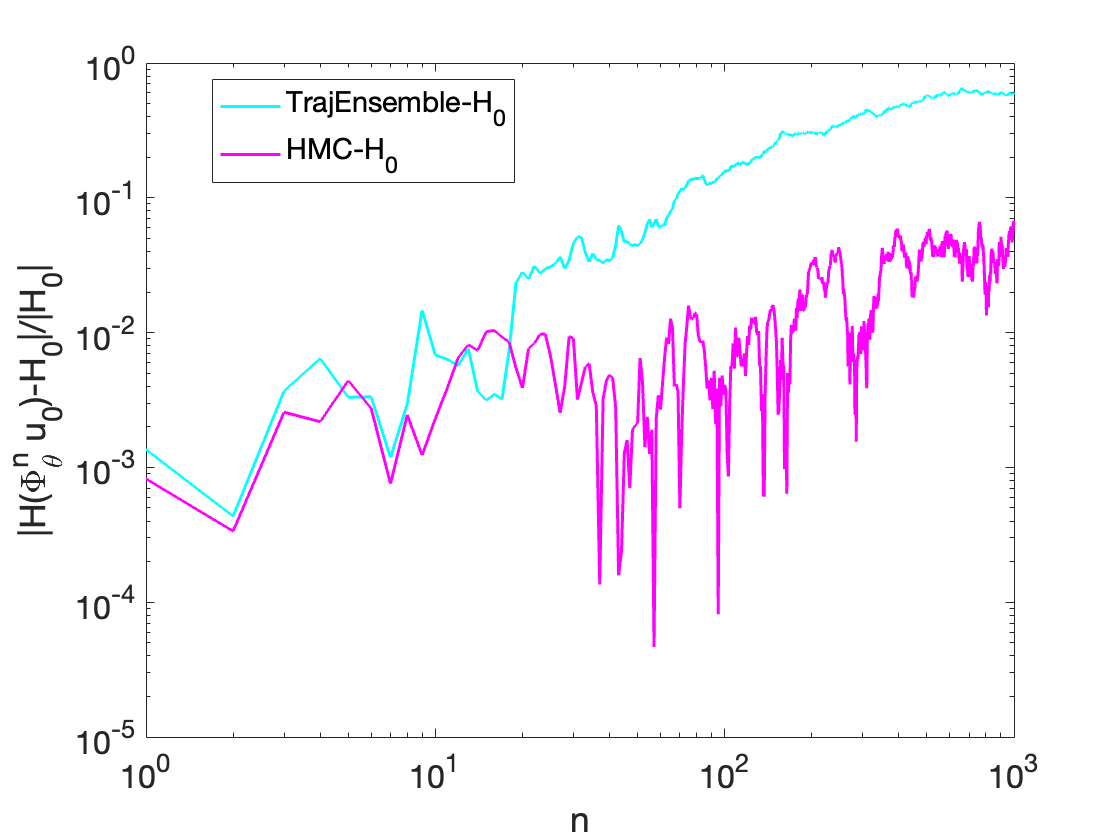}
        \caption{energy error}
    \end{subfigure}
    \caption{(FPUT problem, $\omega=300$) Effect of the sampling mechanism at fixed target energy. (a) Minimum distance from each point of the reference trajectory to the set of training initial conditions, for the HMC-$H_0$ and TrajEnsemble-$H_0$ datasets. (b) Energy errors of the long-time trajectories generated by $\Phi_{T_0=1}$ trained on each dataset. 
    }
    \label{fig:fput_data}
\end{figure}

The comparison is sharper on FPUT at $\omega=300$, where both datasets target the same energy level and differ only in the sampling mechanism. Against HMC-$H_0$ we place TrajEnsemble-$H_0$, the natural alternative: sample momenta, integrate each initial state deterministically, and collect every state along the resulting trajectories, appealing to ergodicity to cover the level set. Figure~\ref{fig:fput_data}~(b) shows that the map trained on HMC-$H_0$ data attains markedly better long-time energy preservation. Figure~\ref{fig:fput_data}~(a) indicates why: the minimum distance from each point of the reference trajectory to the set of training initial conditions stays small and stationary for the HMC-$H_0$ dataset, while for TrajEnsemble-$H_0$ it is larger and grows along the trajectory, so the learned map is evaluated progressively further from the region in which it was fitted. This behavior is consistent with Proposition~\ref{prop:hmc-invariance}: the HMC-$H_0$ transition leaves the target microcanonical distribution invariant, whereas a finite deterministic trajectory ensemble can retain substantial dependence on the states from which its trajectories were launched.

\subsubsection*{Effect of the slow-block Taylor order}
Lemma~\ref{lem:consistency} establishes the built-in consistency of the activated Taylor architecture. We therefore focus the numerical ablation on the slow--fast FPUT construction actually used in the principal multiscale experiments.

FPUT at $\omega=50$ tests the slow--fast variant, in which different orders are assigned to the slow and fast blocks. Here the time collocation points are uniform on $[0,0.125]$, the phase-space points are drawn from the HMC-$H_0$ dataset of Section~\ref{sec:res_phasespace}, the residual uses velocity Verlet with $h=2^{-10}$, and each model is trained with batch size $2000$ for 1M iterations; the final losses are $1.7\times10^{-4}$ (train) and $1.8\times10^{-4}$ (test) for $(p_f,p_s)=(0,0)$, and $1.8\times10^{-4}$ and $1.9\times10^{-4}$ for $(0,2)$.

\begin{figure}[ht]
    \centering
    \includegraphics[width=0.7\linewidth]{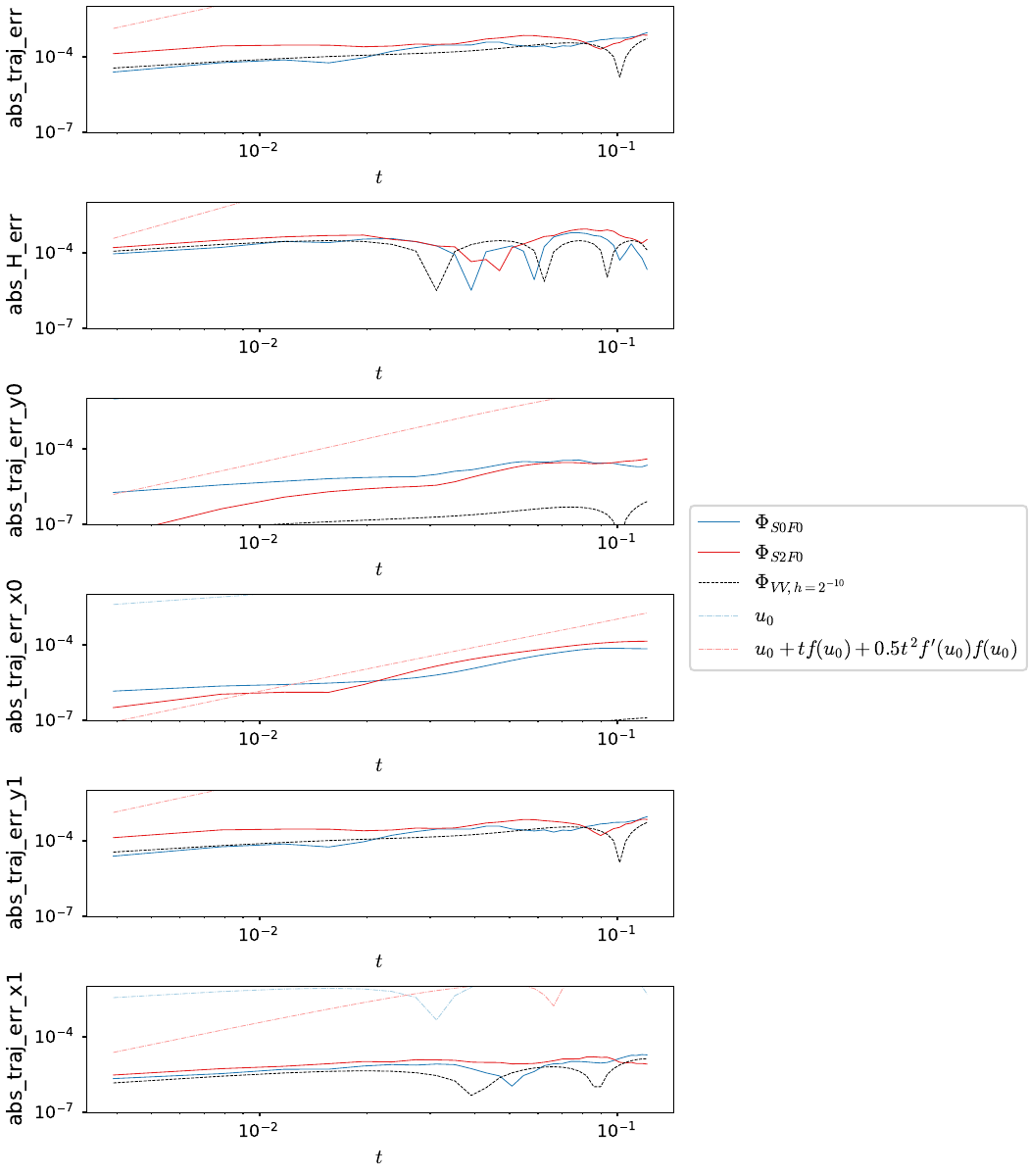}
    \caption{(FPUT problem, $\omega=50$) one-step prediction errors of slow--fast Taylor-based flow maps with truncation orders $(p_f,p_s)=(0,0)$ and $(0,2)$.}
    \label{fig:fput_taylororder_1step}
\end{figure}

\begin{figure}[ht]
    \centering
    \begin{subfigure}[b]{0.35\linewidth}
        \includegraphics[width=\linewidth,height=12cm,keepaspectratio,trim=0 0 110 0, clip]{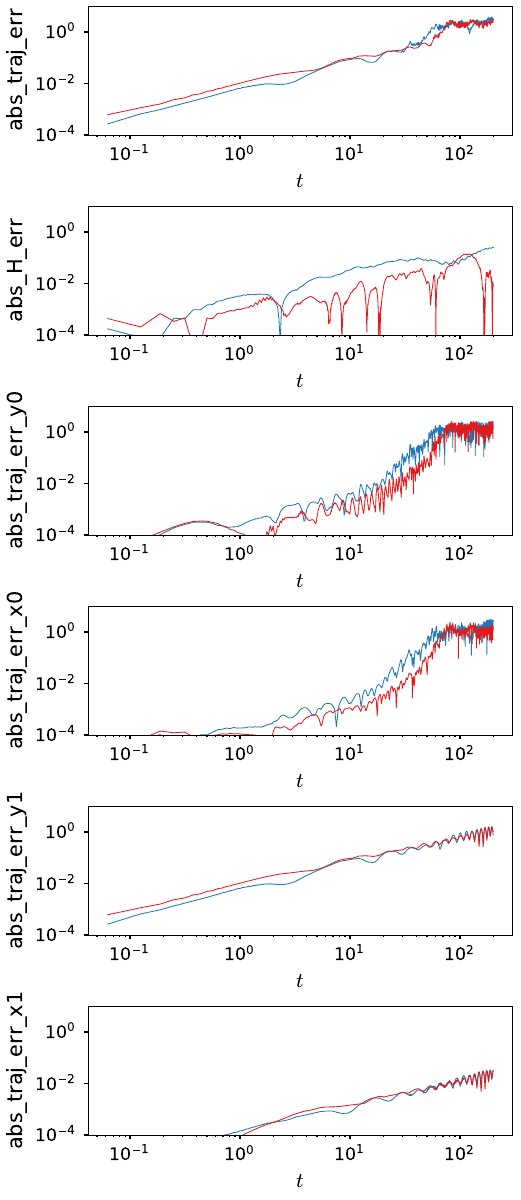}
    \end{subfigure}
    \begin{subfigure}[b]{0.5\linewidth}
        \includegraphics[width=\linewidth,height=12cm,keepaspectratio]{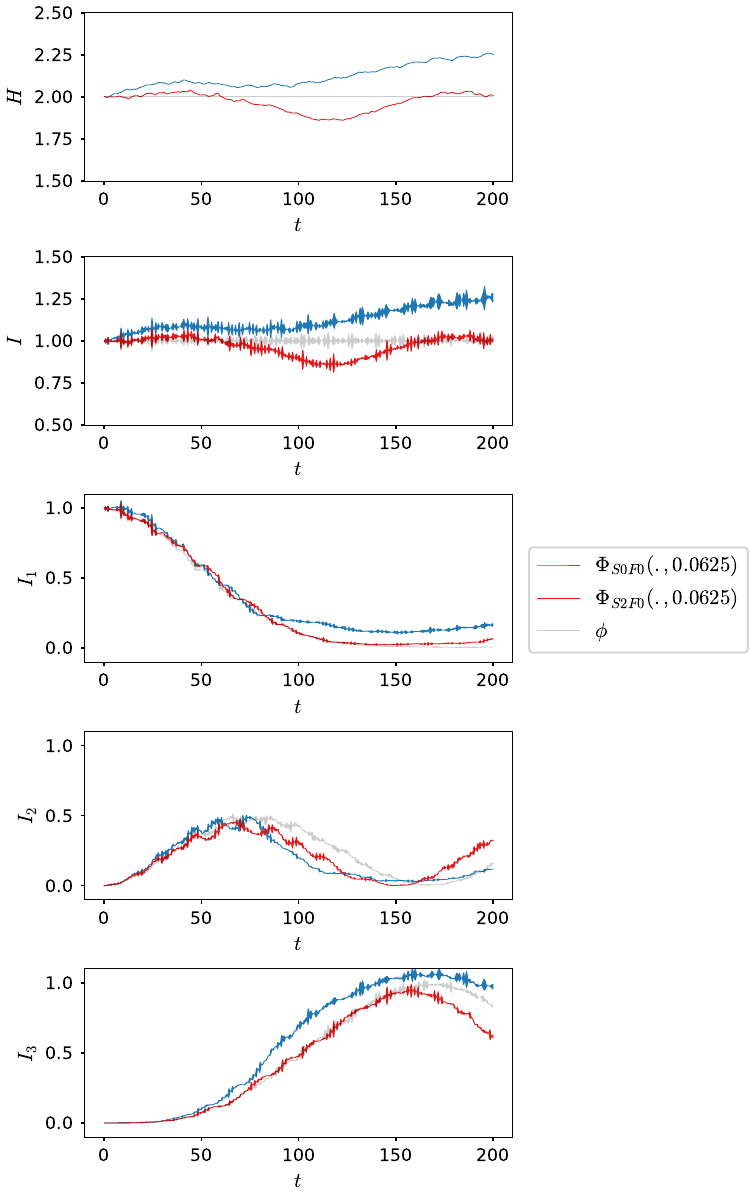}
    \end{subfigure}
    \caption{(FPUT problem, $\omega=50$) Errors and energy profiles of long-time trajectories generated with a fixed step size $t=0.0625$ by slow--fast Taylor-based flow maps with truncation orders $(p_f,p_s)=(0,0)$ and $(0,2)$.}
    \label{fig:fput_taylororder_longtime}
\end{figure}

Figure~\ref{fig:fput_taylororder_1step} shows that the overall one-step errors of the two models are nearly indistinguishable and that the fast components $y_1,x_1$ are captured equally well by both, at the accuracy of the underlying scheme. The difference is confined to the slow components $y_0,x_0$, where $(0,2)$ is clearly more accurate. That difference becomes substantially more consequential under recursion: Figure~\ref{fig:fput_taylororder_longtime} shows long-time trajectories generated with a fixed step $t=0.0625$ up to $T=200$, where the $(0,2)$ model reproduces the energy exchange among the stiff springs and reduces the energy error substantially over $[0,50]$. Thus, in this FPUT experiment, increasing the slow-block order from $p_s=0$ to $p_s=2$ while keeping $p_f=0$ improves the slow components and long-time energy behaviour without materially changing the fast-component accuracy. We did not vary $p_f$, so these runs say nothing about whether a higher order on the fast block would help.

\subsubsection{The training horizon}\label{sec:res_timewindow}
\begin{figure}[ht]
    \centering
    \includegraphics[width=0.7\linewidth]{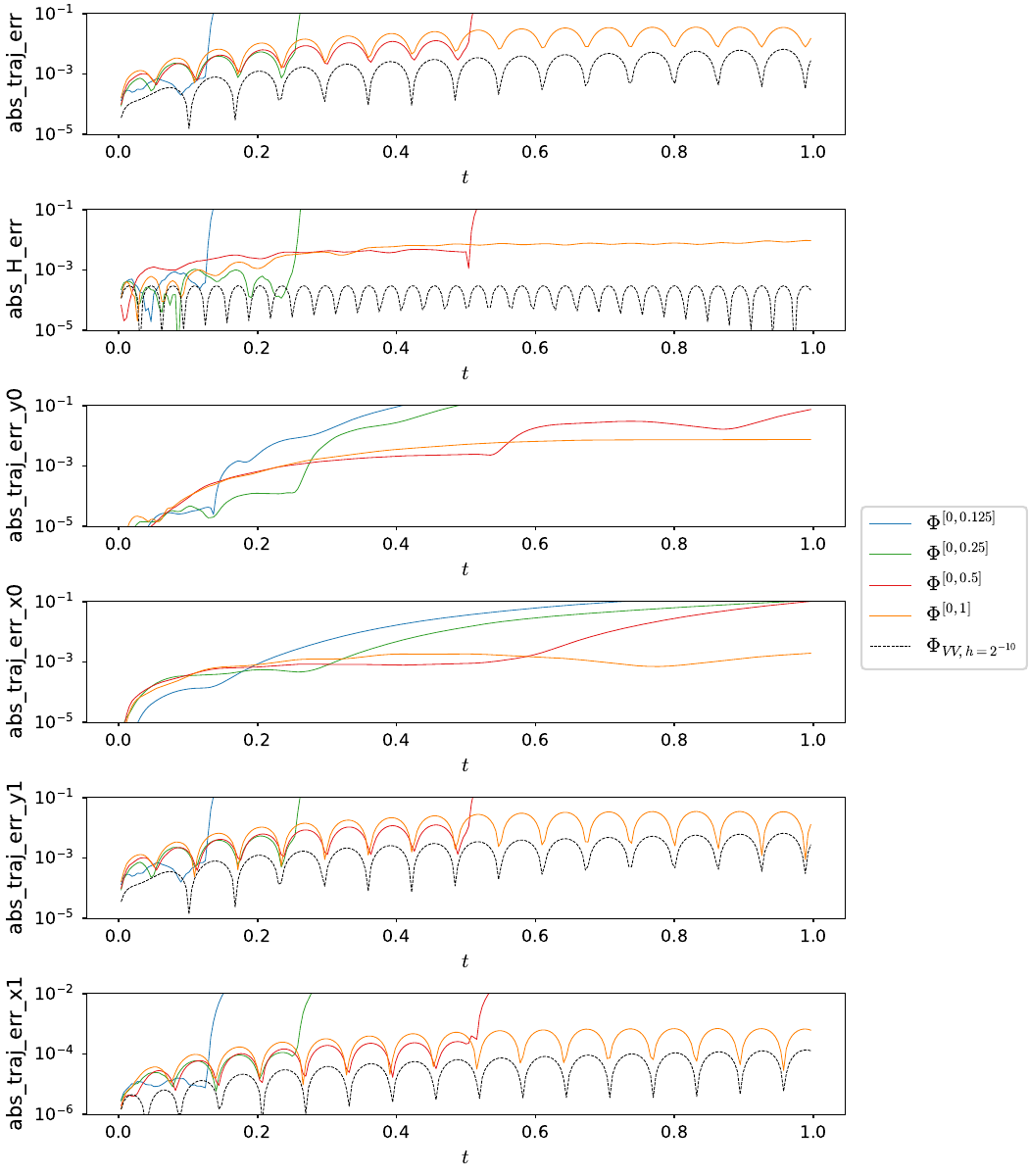}
    \caption{(FPUT problem, $\omega=50$) one-step prediction errors of flow maps learned up to $T=0.125,0.25,0.5,1$.}
    \label{fig:fput_timewindowlength_1step}
\end{figure}

\begin{figure}[ht]
    \centering
    \begin{subfigure}[b]{0.35\linewidth}
        \includegraphics[width=\linewidth,height=12cm,keepaspectratio,trim=0 0 115 0, clip]{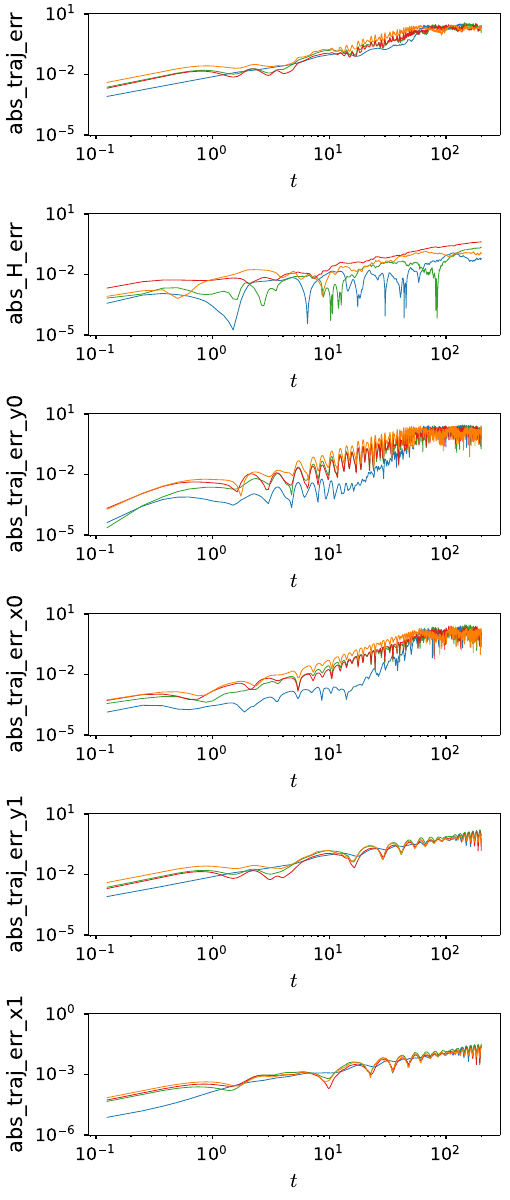}
    \end{subfigure}
    \begin{subfigure}[b]{0.5\linewidth}
        \includegraphics[width=\linewidth,height=12cm,keepaspectratio]{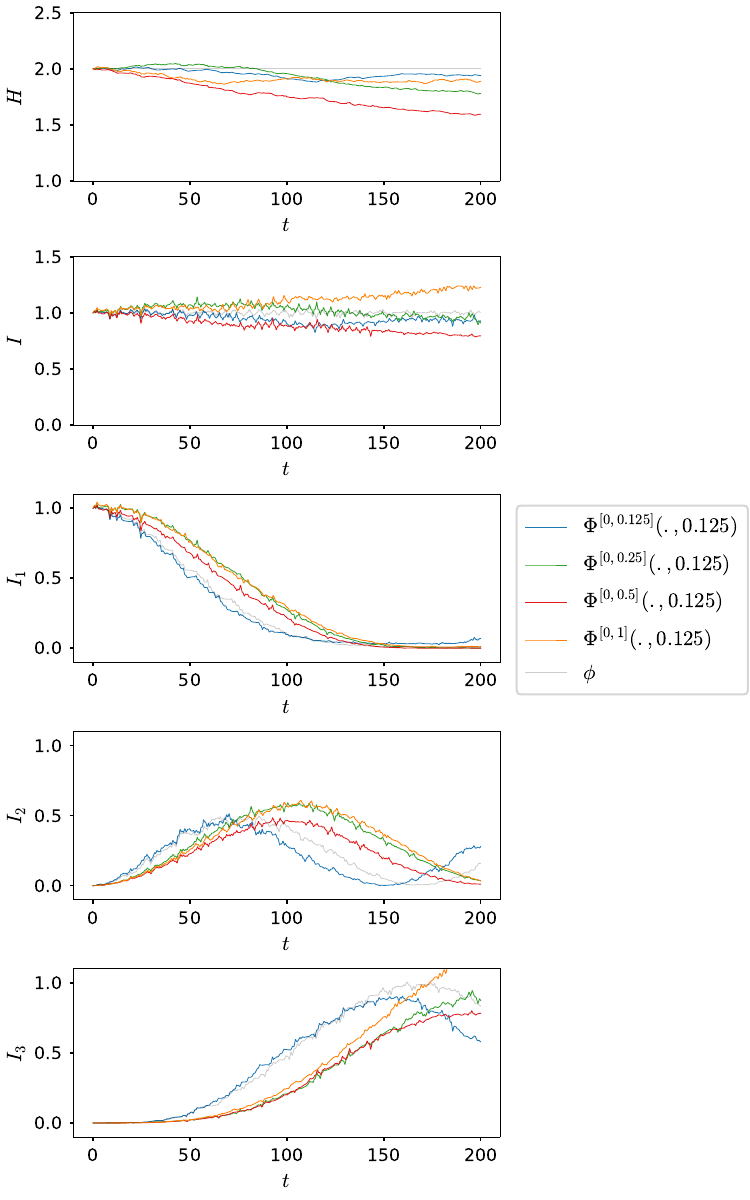}
    \end{subfigure}
    \caption{(FPUT problem, $\omega=50$) Errors and energy profiles of long-time trajectories over $[0,200]$ generated with a fixed step size $t=0.125$ by flow maps learned up to $T=0.125,0.25,0.5,1$.}
    \label{fig:fput_timewindowlength_longtime}
\end{figure}

The terminal time $T$ over which a variable-timestep map is learned cannot be enlarged indefinitely. For FPUT at $\omega=50$, we compare $T=0.125,0.25,0.5,1$ with phase-space points from the HMC-$H_0$ dataset and a velocity Verlet residual with $h=2^{-10}$. Figure~\ref{fig:fput_timewindowlength_1step} shows that the shortest horizons give the most accurate slow components, while the longer horizons extend the reach of a single application at the cost of accuracy everywhere within it. Under recursion with step $t=0.125$ over $[0,200]$ (Figure~\ref{fig:fput_timewindowlength_longtime}) the $T=0.125$ model keeps the slow components accurate for longest and gives the closest match to the reference energy profile. We therefore select the smallest horizon compatible with the intended step size, which is the choice reported in the experiments above.

\subsubsection{The multi-step data loss}\label{sec:res_multistep}
We assess the effect of the sequence length $S$ in \eqref{def:data-loss} on FPUT at $\omega=300$ with $T_0=1$, training three models with $S=1,3,5$. The initial states are drawn from the same HMC-$H_0$ dataset, the targets are generated with a high-order symplectic integrator, and all models share the architecture and training budget. Their final train/test losses are comparable: $2.0\times10^{-7}/2.2\times10^{-7}$ for $S=1$, $1.8\times10^{-7}/1.9\times10^{-7}$ for $S=3$, and $2.8\times10^{-7}/3.0\times10^{-7}$ for $S=5$.

\begin{figure}[ht]
    \centering
    \begin{subfigure}[b]{0.35\linewidth}
        \includegraphics[width=\linewidth,height=12cm,keepaspectratio,trim=0 0 70 0, clip]{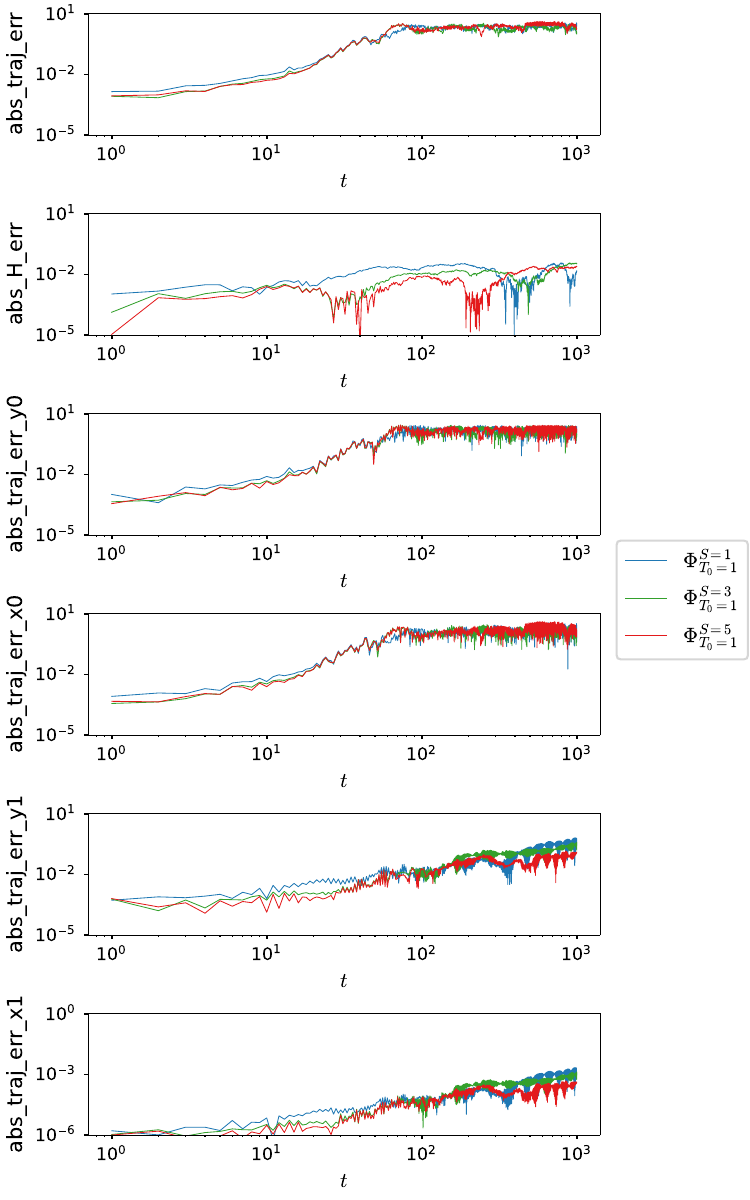}
    \end{subfigure}
    \begin{subfigure}[b]{0.5\linewidth}
        \includegraphics[width=\linewidth,height=12cm,keepaspectratio]{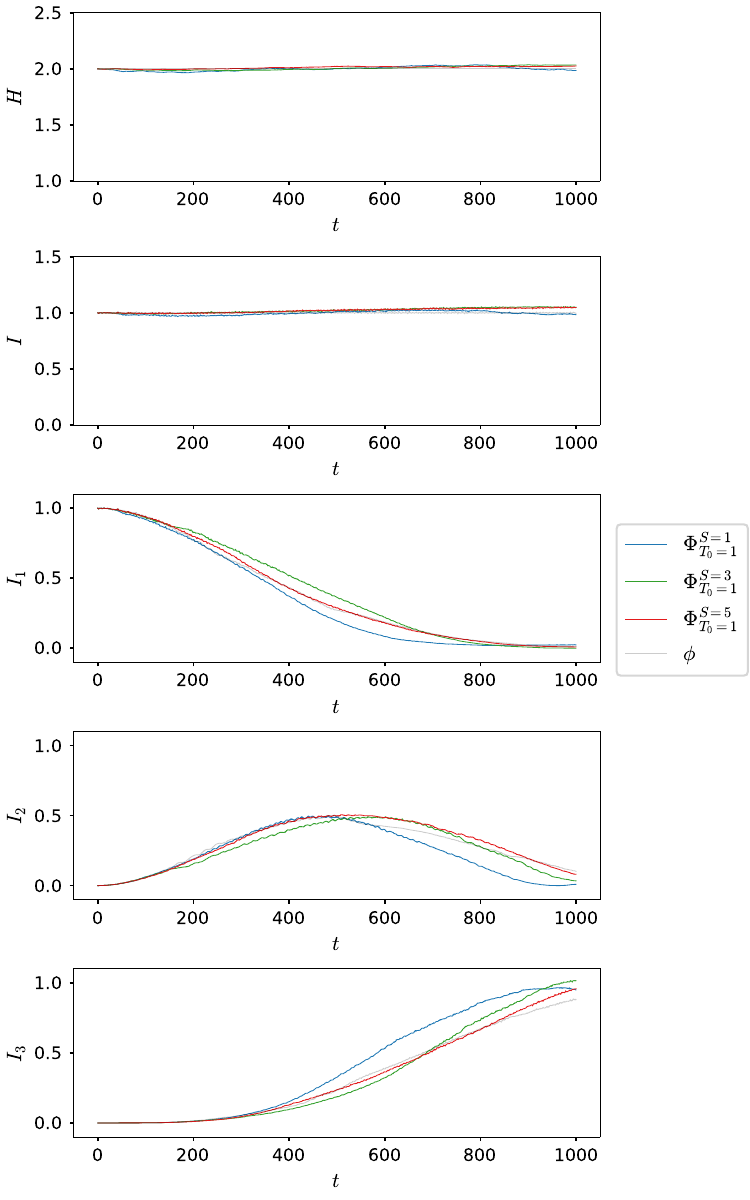}
    \end{subfigure}
    \caption{(FPUT problem, $\omega=300$) Errors and energy profiles of long-time trajectories over $[0,1000]$ generated from an in-distribution initial condition by fixed-timestep flow maps $\Phi_{T_0=1}$ trained with sequence lengths $S=1,3,5$.}
    \label{fig:multistep_id}
\end{figure}

The difference appears under recursion. Figure~\ref{fig:multistep_id} shows that all three have similar trajectory errors up to $T=1000$, while the $S=5$ model gives the lowest energy error over the initial long-time window and best reproduces the reference energy profile. This is the reason the $\omega=300$ experiments of Section~\ref{subsec:FPU} use $S=5$.

\subsection{Positioning with respect to existing learned integrators}\label{sec:baselines}

\paragraph{A parameter-matched symplectic-network comparison}
On the NPCO problem we compare the fixed-timestep map $\Phi^{80\mathrm{k}}_{T_0=5}$ with G-SympNet \cite{jin2020sympnets}, a representative fixed-step symplectic architecture with a publicly available implementation. Both models are trained on the same HMC-$H_0$ dataset, matched at approximately $80$k trainable parameters, and given the same training budget. Apart from its reduced width, $\Phi^{80\mathrm{k}}_{T_0=5}$ uses the same architecture, loss, and data as the 1M-parameter $\Phi_{T_0=5}$ model reported elsewhere.

\begin{figure}[ht]
    \centering
    \begin{subfigure}[b]{0.31\linewidth}
        \includegraphics[width=\linewidth,trim=170 0 0 15, clip]{figures/t_5_phase_space_fixedH_ref.pdf}
        \caption{reference $\phi$}
    \end{subfigure}
    \begin{subfigure}[b]{0.31\linewidth}
        \includegraphics[width=\linewidth,trim=170 0 0 15, clip]{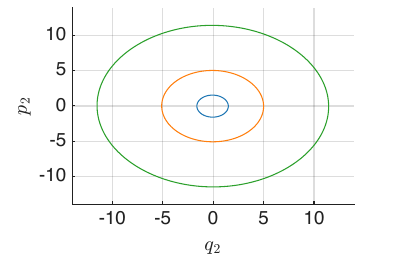}
        \caption{G-SympNet}
    \end{subfigure}
    \begin{subfigure}[b]{0.31\linewidth}
        \includegraphics[width=\linewidth,trim=170 0 0 15, clip]{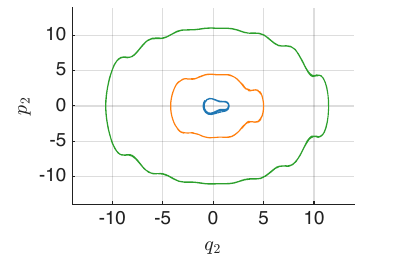}
        \caption{$\Phi^{80\mathrm{k}}_{T_0=5}$}
    \end{subfigure}
    \caption{(NPCO problem) Slow-oscillator phase-space portraits over $[0,2000]$ for the reference flow and the two approximately $80$k-parameter models used in the matched comparison.}
    \label{fig:npco_structure_preserving}
\end{figure}

Figure~\ref{fig:npco_structure_preserving} compares long-time phase portraits over $[0,2000]$. The modified MLP recovers the qualitative structure of the slow-oscillator dynamics more accurately, while G-SympNet preserves the Hamiltonian more faithfully. Table~\ref{tab:baseline-comparison} quantifies the same tradeoff: at matched capacity the modified MLP has smaller one-step and long-time trajectory errors, whereas G-SympNet has smaller long-time energy error. Thus symplectic structure alone does not guarantee smaller trajectory error at a fixed model capacity.

\begin{table}[ht]
    \centering
    \small
    \caption{Parameter-matched comparison on the NPCO problem. Both models use identical HMC-$H_0$ training data and the same training budget.}
    \label{tab:baseline-comparison}
    \begin{tabularx}{\textwidth}{l XX XX X}
        \toprule
        \multirow{2}{*}{\textbf{Model}}
            & \multicolumn{2}{c}{\textbf{one-step}}
            & \multicolumn{2}{c}{\textbf{Long-time} ($T{=}2000$)}
            & \multirow{2}{*}{\textbf{ms/step}} \\
        \cmidrule(lr){2-3} \cmidrule(lr){4-5}
            & \textbf{Traj err} & \textbf{$H$ err}
            & \textbf{Traj err} & \textbf{$H$ err} & \\
        \midrule
        $\Phi^{80\mathrm{k}}_{T_0=5}$ (ours)
            & $8.4\times 10^{-3}$ & $6.7\times 10^{-3}$
            & $2.4$              & $1.2\times 10^{-1}$
            & $0.21$ \\
        G-SympNet
            & $8.8\times 10^{-2}$ & $5.7\times 10^{-3}$
            & $3.0$              & $3.0\times 10^{-3}$
            & $1.3$ \\
        \bottomrule
    \end{tabularx}
\end{table}

\paragraph{Scope of existing learned integrators}
Vector-field learners approximate $H$ or $f$ and still require numerical integration at inference, so their prediction cost remains tied to the chosen integrator; stiff-system variants include mono-implicit Runge--Kutta training on FPUT-type systems \cite{noren2023learning} and multi-resolution approaches \cite{li2026frequency}. Flow-map learners instead advance the state directly over a coarse increment. Published examples of the latter generally involve more moderate scale separation \cite{jin2020sympnets,burby2020fast,duruisseaux2023approximation}; the FPUT experiment at $\omega=300$ tests a coarse time-$1$ map across substantially faster oscillations. For the $\alpha$-particle problem, our setup follows \cite{burby2024non}, but the comparison is methodological rather than head-to-head: that work learns an adiabatic invariant, whereas we learn a single $\epsilon$-parametric map for the full reduced dynamics.

\subsection{Runtime comparison}\label{subsec:runtime-comp}

Table~\ref{tab:runtime_longtime} summarizes the runtime and accuracy tradeoffs between our neural network solvers and conventional numerical solvers.
The neural network solvers are implemented in Python and run on an NVIDIA Tesla V100 GPU. The numerical solvers are implemented in Julia and run on an Intel Xeon Gold 6248R CPU. The batch evaluations for numerical solvers are run in parallel on 16 cores.
For each problem, we focus on a short-time integration interval $T_s$ and compare the trajectory error, energy error, and runtime for both a single trajectory and for long-time integration from a batch of initial conditions (e.g., 512 trajectories up to $1000T_s$).

Throughout, ``runtime'' refers to \emph{inference} wall-clock time only: the time to advance the stated number of trajectories over the stated horizon by repeated application of the trained network or of the numerical integrator. Training costs are excluded; they are reported in Table~\ref{tab:flow_map_training}, and a learned flow map is worthwhile only when its (one-time) training cost can be amortized over many ensemble members, long horizons, or repeated use. We also note that the comparison spans different hardware and software stacks (GPU/Python for the networks, 16-core CPU/Julia for the integrators), reflecting the natural deployment of each method; absolute ratios should be read with this caveat in mind. The reference integrators used here take fixed steps and are themselves straightforward to batch on a GPU, so the reported ratios are not hardware-independent measurements. The platform-independent algorithmic distinction is that the learned map advances by the coarse increment $T_0$, whereas the reference integrator must take sufficiently small resolving steps. The reference integrators (velocity Verlet, implicit midpoint) are the standard symplectic choices for these problems; specialized multiscale integrators (e.g., trigonometric methods for oscillatory problems \cite{geometric}) could narrow the reported gaps for the FPUT problem.

The clearest advantage arises for the FPUT problem with $\omega=300$, where $T_s$ is large compared to the smallest timescale $\omega^{-1}$: the fixed-timestep flow map $\Phi_{T_0=1}$ advances 512 trajectories to $1000\,T_s$ in $1.4$\,s, versus $265$\,s for velocity Verlet at the step size ($h=2^{-15}$) needed to resolve the fast oscillations---a ${\sim}190\times$ speedup at comparable trajectory and energy errors. For a fixed trained architecture, the cost of one network evaluation does not explicitly depend on $\omega$, whereas the number of small steps required by the resolving reference integrator grows with the scale separation; this makes the learned-map advantage potentially larger as the stiffness increases. The advantage is also largest in the ensemble setting, where GPU batching efficiently amortizes the cost of network evaluation across the ensemble. By contrast, for the NPCO and $\alpha$-particle problems, where $T_s$ is comparable to the fast timescale and the dimension is low, well-implemented conventional integrators remain faster for single trajectories, and the learned maps are only competitive at the ensemble level. We report these unfavorable cases deliberately: they help delineate regimes in which learned flow maps can be advantageous---strong timescale separation, large ensembles, or repeated simulation---rather than serving as a universal replacement for numerical integrators.

\begin{table}[ht]
\centering
\footnotesize
\caption{Runtime and accuracy comparison of various solvers for short-time integration from a single initial condition and long-time integration from a batch of initial conditions. For the $\alpha$-particle problem, the learned flow map preserves $\|v\|$---and hence the conserved energy---exactly by construction (Section~\ref{sec:runtime}), so its energy error is identically zero.
}
\label{tab:runtime_longtime}
\setlength{\tabcolsep}{4pt}
\begin{tabular}{lP{2.2cm}C{0.5cm}C{1.5cm}C{1.7cm}C{1.4cm}C{1.9cm}}
\toprule
\multirow{2}{*}{\textbf{Problem}} & \multirow{2}{*}{\textbf{Solver}} & \multirow{2}{*}{$\bm{T_s}$}
& \multicolumn{3}{l}{\textbf{1 Trajectory to \( T_s \)}} 
& \textbf{512 Traj.'s to  $1000 T_s$ } \\
\cmidrule{4-7}
& & & Traj err  & $H$ error & Runtime & Runtime \\
\midrule
\multirow{2}{*}{NPCO} 
  & $\Phi(.,5)$              & 20 & $5.3\times10^{-4}$ & $1.7\times10^{-5}$ & 13 ms & 14 s \\
  & $\Phi_{VV,h=0.02}$       & 20 & $4.9\times10^{-4}$ & $1.0\times10^{-4}$ & 0.35 ms & 12 s\\
\midrule
\multirow{3}{*}{FPUT $\omega=50$}
  & $\Phi(.,2^{-3})$         & 1 & $7.8\times10^{-3}$  & $2.7\times10^{-4}$ & 33 ms & 35 s\\
  & $\Phi_{T_0=1}$           & 1 & $1.7\times10^{-3}$ & $4.5\times10^{-4}$ & 1.2 ms & 1.2 s \\
  & $\Phi_{VV,h=2^{-11}}$    & 1 & $9.8\times10^{-4}$ & $3.0\times10^{-5}$ & 0.53 ms & 17 s\\
\midrule
\multirow{2}{*}{FPUT $\omega=300$} 
  & $\Phi_{T_0=1}$           & 1 & $8.7\times10^{-4}$ & $4.5\times10^{-8}$ & 1.4 ms & 1.4 s\\
  & $\Phi_{VV,h=2^{-15}}$    & 1 & $1.1\times10^{-3}$ & $2.3\times10^{-7}$ & 8.3 ms & 265 s\\
\midrule
\multirow{2}{*}{$\alpha$-particle} 
  & $\Phi(.,5;0.27)$         & 10 & $6.8\times10^{-4}$ & 0 & 4.6 ms & 4.9 s \\
  & $\Phi_{Midpt,h=2^{-9}}$  & 10 & $5.1\times10^{-4}$ & $1.3\times10^{-4}$ & 0.16 ms & 5.1 s \\
\bottomrule
\end{tabular}
\end{table}

\section{Conclusion}\label{sec:conc}

We developed a framework for learning fixed- and variable-timestep Hamiltonian flow maps using either numerical-scheme residuals or reference trajectory data. The main components are Taylor-consistent neural parameterizations, energy-balanced losses for stiff variables, and microcanonical HMC-$H_0$ sampling. The critical-point analysis shows precisely what residual minimization enforces at the function-space level: under local invertibility of the implicit update, every $C^1$ critical point of the unrestricted residual functional has zero residual and therefore reproduces the chosen numerical scheme. The multiscale error analysis shows that, after the energy-based rescaling, stiff linear oscillations do not enter the error-amplification rate; linear error accumulation follows while the slow-scale amplification remains small.

The experiments indicate two useful operating regimes. For a prescribed coarse increment $T_0$, a fixed-timestep map trained from accurate reference data is the simplest and often the most accurate surrogate. When predictions are required at arbitrary intermediate times, the Taylor-consistent variable-timestep map can instead be trained from a scheme residual, using a method appropriate to the Hamiltonian structure. For slow--fast systems, raising the truncation order on the slow block improved both the slow components and the long-time energy behavior in our FPUT experiment, with the fast block held at order zero; higher orders on the fast block were not tested. For the oscillatory class considered here, energy-balanced scaling prevents errors in fast positions from being underweighted; in the FPUT ablation, this weighting is essential for resolving the fast energy exchange. HMC-$H_0$ is most useful when training should concentrate on a target energy band. The runtime comparisons also show that the learned map is not uniformly preferable to a conventional integrator: its advantage appears when the cost of resolving fast scales is high and the trained surrogate is reused across sufficiently large trajectory ensembles.

The principal limitations are inherited from the training target and from dynamical sensitivity. A residual-trained map approximates the selected numerical scheme, so its discretization error and stability properties remain relevant. The stiffness-uniform theorem does not remove exponential sensitivity associated with slow nonlinear or chaotic dynamics, and the current analysis does not cover systems such as the $\alpha$-particle model whose fast linearization varies strongly with position. Extending the error analysis to such noncanonical systems, and combining the present loss and sampling ideas with explicitly symplectic architectures, are natural directions for future work. The present benchmarks are also moderate-dimensional, and the scaling of the approach with phase-space dimension remains to be investigated.

\section*{Code and reproducibility}
The implementation accompanying this paper is publicly available at \url{https://github.com/fr0420/learnsolnmap}.

\section*{Acknowledgments}
The authors thank Bjorn Engquist and  Joshua Burby for motivating conversations.
The authors are partially supported by NSF Grant DMS-2208504. Tsai is partially supported by NSF Grant DMS-2513857.

\bibliographystyle{abbrv} 
\bibliography{references}

\appendix

\section{Critical points of the residual functional}\label{sec:proofs_critical_points}

We prove Theorem~\ref{thm:general-one-step}. The argument is continuous in time and avoids decomposing the integral into shifted discrete grids. The time density $\gamma$ is carried throughout; setting $\gamma\equiv1$ gives the unweighted functional.

Fix a critical point $\Phi$ and define, for each $u$,
\[
    w(u,t):=Q\lr{\Phi(u,t+h),\Phi(u,t)},
\]
\[
    A(u,t):=\partial_1Q\lr{\Phi(u,t+h),\Phi(u,t)},
    \quad
    B(u,t):=\partial_2Q\lr{\Phi(u,t+h),\Phi(u,t)}.
\]
For an admissible $C^1$ variation $\psi$ with $\psi(u,0)=0$, the first variation is
\begin{align}
\delta\mathcal R[\Phi](\psi)
  &=\int_\Omega\int_0^T\gamma(t)
     \innerp{w(u,t)}{A(u,t)\psi(u,t+h)+B(u,t)\psi(u,t)}
     \,dt\,d\rho(u)\nonumber\\
  &=\int_\Omega\biggl[
     \int_h^{T+h}\gamma(s-h)\innerp{A(u,s-h)^Tw(u,s-h)}{\psi(u,s)}\,ds\nonumber\\
  &\hspace{4em}
     +\int_0^T\gamma(s)\innerp{B(u,s)^Tw(u,s)}{\psi(u,s)}\,ds
     \biggr]\,d\rho(u),
\label{eq:continuous-first-variation}
\end{align}
where the second identity uses $s=t+h$ in the first integral.

Because $\rho$ has full support and all coefficients are continuous, the fundamental lemma of the calculus of variations can be applied pointwise in $u$. Variations supported in $s\in(T,T+h)$ meet the first integral only, and give
\begin{equation}\label{eq:continuous-terminal}
    \gamma(s-h)\,A(u,s-h)^Tw(u,s-h)=0,
    \qquad T<s<T+h.
\end{equation}
Since $\gamma>0$ and $A$ is invertible, $w(u,t)=0$ for $t\in(T-h,T)$; here $T\ge h$ guarantees that $(T,T+h)$ lies in the range $(h,T+h)$ of the shifted integral.

For $h<s<T$, variations supported near $s$ give the interior adjoint relation
\begin{equation}\label{eq:continuous-adjoint}
    \gamma(s-h)A(u,s-h)^Tw(u,s-h)+\gamma(s)B(u,s)^Tw(u,s)=0.
\end{equation}
Starting from the terminal interval on which $w=0$, equation \eqref{eq:continuous-adjoint} propagates zero backward by intervals of length $h$: if $w(u,s)=0$ on an interval, the second term vanishes there, and $\gamma(s-h)>0$ together with the invertibility of $A(u,s-h)$ forces $w(u,s-h)=0$ on the preceding interval. Repeating this argument reaches every $t\in(0,T)$. Continuity of $w$ then gives $w(u,0)=0$ as well. Hence $R_h[\Phi](u,t)=0$ on $\Omega\times[0,T]$, proving Theorem~\ref{thm:general-one-step}.

Positivity of $\gamma$ enters only through these two steps, and it is exactly what fails for a time measure concentrated on a discrete set.

\section{Proof of the stiffness-uniform error bound}
\label{sec:proof_global_error}

Throughout this appendix, $u=(p,q)\in\mathbb{R}^{2d}$, and the vector field of the Hamiltonian \eqref{eq:osc-class} is
$f(u)=\lr{-\Omega^2q-\nabla U(q),\,p}$.
As in Section~\ref{sec:ana_errorbound}, let
$\tilde\Omega=\operatorname{diag}(\tilde\omega_1,\ldots,\tilde\omega_d)$ with
$\tilde\omega_j=\max\{\omega_j,1\}$.

\subsection{Proof of \texorpdfstring{Theorem~\ref{thm:linear-error-growth}}{the stiffness-uniform theorem}}
\label{sec:proof-linear-growth-main}

\begin{proof}
Set $e_n=\tilde u_n-u_n$. Since
$\tilde u_{n+1}=\phi_{\Delta t}(\tilde u_n)+\delta(\tilde u_n)$,
\begin{align}\label{eq:error-recurrence-integral}
 e_{n+1}
 =\lr{\int_0^1 D\phi_{\Delta t}\lr{u_n+\theta e_n}\,d\theta}e_n
 +\delta(\tilde u_n).
\end{align}
Thus it remains to bound the derivative of the exact flow in the weighted norm.

\emph{Step 1: rescaling the variational equation.}
For \eqref{eq:osc-class},
\[
Df(u)=
\begin{pmatrix}0&-\Omega^2-\nabla^2U(q)\\ I&0\end{pmatrix}
=A+P(u),
\]
where
\begin{align}\label{eq:A-P-decomposition}
A&=\begin{pmatrix}0&-\tilde\Omega^2\\ I&0\end{pmatrix},
&
P(u)&=\begin{pmatrix}0&B(q)\\0&0\end{pmatrix},
&
B(q)&=\tilde\Omega^2-\Omega^2-\nabla^2U(q).
\end{align}
The correction $\tilde\Omega^2-\Omega^2$ is supported only on coordinates with $\omega_j<1$. Hence Assumption~\ref{ass:system-class} gives
\begin{equation}\label{eq:B-bound}
 \|B(q)\|\le 1+c_U=:M
 \qquad\text{on }\mathcal E'_{H_0}.
\end{equation}
Define
\[
S=\operatorname{diag}(I,\tilde\Omega^{-1}),
\qquad \|S^{-1}v\|_2=\|v\|_W.
\]
Then
\[
S^{-1}AS=
\begin{pmatrix}0&-\tilde\Omega\\ \tilde\Omega&0\end{pmatrix},
\]
so
\[
R(t):=S^{-1}e^{At}S=
\begin{pmatrix}
\cos(\tilde\Omega t)&-\sin(\tilde\Omega t)\\
\sin(\tilde\Omega t)&\cos(\tilde\Omega t)
\end{pmatrix}
\]
is orthogonal. Moreover,
\begin{equation}\label{eq:Phat-bound}
\hat P(u):=S^{-1}P(u)S
=\begin{pmatrix}0&B(q)\tilde\Omega^{-1}\\0&0\end{pmatrix},
\qquad
\|\hat P(u)\|\le M,
\end{equation}
with $M$ independent of the stiff frequencies.

We also record a uniform energy-gradient bound. On $\mathcal E'_{H_0}$,
\begin{equation}\label{eq:scaled-energy-gradient}
\|S^T\nabla H(u)\|_2
=\norm{\lr{p,\tilde\Omega^{-1}(\Omega^2q+\nabla U(q))}}
\le K_H,
\end{equation}
where $K_H$ is independent of $\Omega$; indeed,
$\|\tilde\Omega^{-1}\Omega^2q\|\le\|\Omega q\|$ and
$\|\tilde\Omega^{-1}\nabla U\|\le\|\nabla U\|$.
Choose $E_0>0$ so that $K_HE_0<\Delta H$.
If $H(u)=H_0$ and $\|v-u\|_W\le E_0$, the segment joining $u$ and $v$ stays in $\mathcal E'_{H_0}$: otherwise its first exit would require an energy change $2\Delta H$, contradicting \eqref{eq:scaled-energy-gradient} and the mean-value theorem. In particular, every point on that segment has energy within $\Delta H$ of $H_0$, and its exact trajectory remains in $\mathcal E'_{H_0}$ by energy conservation.

\emph{Step 2: stiffness-uniform one-step amplification.}
For a base point $\xi$, let $X(t)=D\phi_t(\xi)$ and $\hat X(t)=S^{-1}X(t)S$. The variation-of-constants formula gives
\begin{equation}\label{eq:variationofconstants}
\hat X(t)=R(t)+\int_0^t R(t-s)\hat P(\phi_s(\xi))\hat X(s)\,ds.
\end{equation}
Whenever $\xi$ lies in the weighted $E_0$-neighborhood above, its exact trajectory remains in $\mathcal E'_{H_0}$. Since $R(t)$ is orthogonal, Gr\"onwall's inequality and \eqref{eq:Phat-bound} yield
\begin{equation}\label{eq:Xhat-gronwall}
 \|\hat X(\Delta t)\|\le e^{M\Delta t}=:e^C,
 \qquad C=M\Delta t.
\end{equation}
Thus, as long as $\|e_n\|_W\le E_0$, \eqref{eq:error-recurrence-integral} and Assumption~\ref{ass:local-error} imply
\begin{equation}\label{eq:hat-recurrence}
 \|e_{n+1}\|_W\le e^C\|e_n\|_W+\delta^*.
\end{equation}

\emph{Step 3: accumulation and energy error.}
Since $e_0=0$, iteration of \eqref{eq:hat-recurrence} gives
\begin{equation}\label{eq:unrolled-recurrence}
 \|e_n\|_W
 \le \frac{e^{nC}-1}{e^C-1}\,\delta^*
 =G_n(C)\delta^*.
\end{equation}
If $G_N(C)\delta^*\le E_0$, this estimate closes by induction for every $n\le N$. The position component satisfies
\[
 |e_{n,q_j}|\le \tilde\omega_j^{-1}G_n(C)\delta^*.
\]
Using \eqref{eq:scaled-energy-gradient} along the segment from $u_n$ to $\tilde u_n$ gives
\[
 |H(\tilde u_n)-H(u_n)|
 \le K_H\|e_n\|_W
 \le K_HG_n(C)\delta^*,
\]
so the theorem holds with $C_H=K_H$. Finally, if $NC\le c_0$, then
\[
G_n(C)=\sum_{k=0}^{n-1}e^{kC}\le e^{c_0}n,
\]
which gives the stated linear-in-$n$ consequence.
\end{proof}

\subsection{Verification of the hypotheses for the benchmark problems}
\label{sec:hypotheses-verification}

\paragraph{FPUT with $m$ stiff springs}
In the slow--fast variables of Section~\ref{subsec:FPU}, the FPUT Hamiltonian \eqref{eq:FPUT-H(y,x)} has the form \eqref{eq:osc-class} after identifying $p=y$ and $q=x$. There are $m$ slow frequencies $\omega_j=0$ and $m$ fast frequencies $\omega_j=\omega=\epsilon^{-1}$, while the nonlinear potential is
\begin{align*}
U(x)=\;&\frac14(x_{s,1}-x_{f,1})^4+\frac14(x_{s,m}+x_{f,m})^4\\
&+\frac14\sum_{i=1}^{m-1}(x_{s,i+1}-x_{f,i+1}-x_{s,i}-x_{f,i})^4.
\end{align*}
Write the linear forms in these quartic terms as $\ell_i(x)=a_i^Tx$, with $\|a_i\|\le2$. Since $U\ge0$ and $\frac14\ell_i^4\le H_0+\Delta H$ on the energy band, the quantities $\ell_i$, and hence $\nabla U=\sum_i\ell_i^3a_i$ and $\nabla^2U=\sum_i3\ell_i^2a_ia_i^T$, are bounded by constants depending on $H_0$, $\Delta H$, and $m$, but not on $\omega$. The same energy bound gives $\|y\|+\|\Omega x\|=\mathcal O(1)$. Thus Assumption~\ref{ass:system-class} holds, and \eqref{eq:W-norm-theory} coincides with the energy-balanced norm \eqref{def:weighted-l2-norm} used for FPUT.

\paragraph{NPCO}
For \eqref{eq:H_npco}, the linear part of the Jacobian is a direct sum of rotation generators with frequencies $1$ and $\epsilon$, so no rescaling is needed. On an energy band, $|q_1|,|p_1|=\mathcal O(1)$ and $|q_2|,|p_2|=\mathcal O(\epsilon^{-1/2})$. Since the entries of $\nabla^2U$ for $U=q_1q_2\sin(2q_1+2q_2)$ are bounded by a constant multiple of $(1+|q_1|)(1+|q_2|)$, the nonlinear part of the Jacobian satisfies $\epsilon\|\nabla^2U\|=\mathcal O(\epsilon^{1/2})$. The same propagator argument therefore applies with a constant uniform in the scale separation.

\section{Training configurations}\label{sec:training-configurations}

The following tables collect the principal architecture, loss, sampling, and optimization choices used for the reported models.

\begin{table}[ht]
    \centering
    \small
    \caption{Sampling datasets used in the training-distribution ablations.}
    \label{tab:training_data_ablation}
    \begin{tabularx}{\textwidth}{l l >{\raggedright\arraybackslash}X}
        \toprule
        \textbf{Problem} & \textbf{Dataset} & \textbf{Parameters} \\
        \midrule
        \multirow{2}{*}{NPCO}
        & $\mathcal{D}_0^{\text{HMC-$H_0$}}$ & $q_0=(0,0)$, $H_0=1.13$, $M=\mathrm{diag}(1,\epsilon^{-1})$, $\sigma=0.1$, $N_{\mathrm{levelsets}}=1000$, $N=600$, $\lambda=64$, $\phi=\Phi_{\text{KahanLi8},h=2^{-4}}$ \\
        & $\mathcal{D}_0^{\text{BoundingBox}}$ & $N_{\text{samples}}=600000$; $q_1,p_1\in[-1.6,1.6]$, $q_2,p_2\in[-16,16]$ \\
        \midrule
        \multirow{2}{*}{FPUT ($\omega=300$)}
        & $\mathcal{D}_0^{\text{HMC-$H_0$}}$ & $q_0=q_{\text{init}}$, $\sigma=0.1$, $N_{\mathrm{levelsets}}=100$, $N=2000$, $\lambda=0.5$, $\phi=\Phi_{\text{CSS4},h=5^{-6}}$ \\
        & $\mathcal{D}_0^{\text{TrajEnsemble-$H_0$}}$ & $q_0=q_{\text{init}}$, $\sigma=0.1$, $N_{\mathrm{levelsets}}=400$, $N_{\mathrm{traj}}=10$, $L=50$, $\delta t=0.1$, $\phi=\Phi_{\text{CSS4},h=5^{-6}}$ \\
        \bottomrule
    \end{tabularx}
\end{table}

\begin{table}[ht]
    \centering
    \small
    \caption{Summary of flow map architectures. The NPCO entry $\Phi_{T_0=5}$ is the 1M-parameter model; the downscaled $80$k variant $\Phi^{80\mathrm{k}}_{T_0=5}$ of Table~\ref{tab:baseline-comparison} differs from it only in width.}
    \label{tab:flow_map_architectures}
    \renewcommand{\arraystretch}{1.2}
    \begin{tabularx}{\textwidth}{l l X l}
        \toprule
        \textbf{Problem} & \textbf{Flow map} & \textbf{Taylor order} & \textbf{Params} \\
        \midrule
        \multirow{2}{*}{NPCOs}
        & $\Phi$ & $p=2$ & 6M \\
        & $\Phi_{T_0=5}$ & -- & 1M \\
        \midrule
        \multirow{3}{*}{FPUT}
        & $\Phi$ ($\omega=50$) & $p_s=2,p_f=0$ & 7.8M \\
        & $\Phi_{T_0=1}$ ($\omega=50$) & -- & 1.8M \\
        & $\Phi_{T_0=1}$ ($\omega=300$) & -- & 1.8M \\
        \midrule
        $\alpha$-particle
        & $\Phi$ & $p=2$ & 7M \\
        \bottomrule
    \end{tabularx}
\end{table}

\begin{table}[ht]
    \centering
    \small
    \caption{Summary of loss functional definitions. The stated loss is the \emph{complete} training objective for each model: models trained with a residual loss do not include the data loss \eqref{def:data-loss}, and vice versa.}
    \label{tab:flow_map_loss_functionals}
    \renewcommand{\arraystretch}{1.2}
    \begin{tabularx}{\textwidth}{l l l l X}
        \toprule
        \textbf{Problem} & \textbf{Flow map} & \textbf{Phase} & \textbf{Time} & \textbf{Loss} \\
        \midrule
        NPCOs
        & $\Phi$ & Box & $\mathcal{U}(0,5)$
        & MSE-residual loss with \\
        & & & & $\Phi^{VV}_h$ (h=0.01) \\
        & $\Phi_{T_0=5}$ & HMC-$H_0$ & $\{5\}$ & MSE-data loss ($S=1$) \\
        \midrule
        \multirow{4}{*}{FPUT} 
        & $\Phi$ ($\omega=50$) & HMC-$H_0$ & $\mathcal{U}(0,2^{-3})$ 
        & EB-residual loss, with \\ 
        & & & & $\Phi^{VV}_h$(h=$2^{-10}$) \\
        & $\Phi_{T_0=1}$ ($\omega=50$) & HMC-$H_0$ & $\{1\}$ 
        & EB-data loss (S=1) \\
        & $\Phi_{T_0=1}$ ($\omega=300$) & HMC-$H_0$ & $\{1\}$ 
        & EB-data loss (S=5) \\
        \midrule
        $\alpha$-particle
        & $\Phi$ & Shell+Box & $\mathcal{U}(0,5)$ 
        & EB-residual loss, with \\ 
        & & & & $\Phi^{IM}_h  (h=0.01)$ \\
        \bottomrule
    \end{tabularx}
\end{table}

\begin{table}[ht]
    \centering
    \small
    \caption{Summary of training setup and results.}
    \label{tab:flow_map_training}
    \renewcommand{\arraystretch}{1.2}
    \setlength{\tabcolsep}{4pt}
    \begin{tabularx}{\textwidth}{l l X l l l l}
        \toprule
        \textbf{Problem} & \textbf{Flow map} & \textbf{Data} & \textbf{Batch} & \textbf{Iter.} & \textbf{Train} & \textbf{Test} \\
        \midrule
        \multirow{2}{*}{NPCOs}
        & $\Phi$ & -- & 1600 & 5M & $1.5\times10^{-7}$ & $1.5\times10^{-7}$ \\
        & $\Phi_{T_0=5}$ & 480k & {500} & {300k} & {$1.3\times10^{-6}$} & {$2.4\times10^{-6}$} \\
        \midrule
        \multirow{3}{*}{FPUT} 
        & $\Phi$ ($\omega=50$) & -- & 2000 & 1M & $1.8\times10^{-4}$ & $1.9\times10^{-4}$ \\
        & $\Phi_{T_0=1}$ ($\omega=50$) & 160k & 400 & 1M & $3.4\times10^{-7}$ & $4.1\times10^{-7}$ \\
        & $\Phi_{T_0=1}$ ($\omega=300$) & 160k & 400 & 3M & $2.8\times10^{-7}$ & $3.0\times10^{-7}$ \\
        \midrule
        $\alpha$-particle
        & $\Phi$ & -- & 400 & 5M & $1.8\times10^{-7}$ & $1.8\times10^{-7}$ \\
        \bottomrule
    \end{tabularx}
\end{table}

\end{document}